\documentclass{article}
\usepackage{amssymb}

\input{tcilatex}
\begin{document}

\begin{center}
\textbf{The regularity properties of nonlocal abstract wave equations }

\textbf{Veli\ B. Shakhmurov}

Department of Mechanical Engineering, Istanbul Okan University, Akfirat,
Tuzla 34959 Istanbul, E-mail: veli.sahmurov@okan.edu.tr,

Institute of Mathematics and Mechanics, Azerbaijan National Academy of
Sciences, AZ1141, Baku, F. Agaev 9

E-mail: veli.sahmurov@gmail.com

\textbf{Rishad Shahmurov}

shahmurov@hotmail.com

University of Alabama Tuscaloosa USA, AL 35487

\textbf{Abstract}
\end{center}

In this paper, the Cauchy problem for linear and nonlinear nonlocal wave
equations are studied.The equation involves a convolution integral operators
with a general kernel operator functions whose Fourier transform are
operator functions defined in Hilbert space $H$ together with some growth
conditions. We establish local and global existence and uniqueness of
solutions assuming enough smoothness on the initial data and the operator
functions$.$ By selecting the space $H$ and the operators, the wide class of
wave equations in the field of physics are obtained.

\textbf{Key Word:}$\mathbb{\ }$Nonlocal equations, Boussinesq equations%
\textbf{,} wave equations, abstract differential equations, Fourier
multipliers

\begin{center}
\bigskip\ \ \textbf{AMS: 35Lxx, 35Qxx, 47D}

\textbf{1}. \textbf{Introduction}
\end{center}

The aim here, is to study the existence and uniqueness of solution of the
initial value problem (IVP) for nonlocal nonlinear abstract wave equat\i on
(WE) 
\begin{equation}
u_{tt}-a\ast \Delta u+A\ast u=\Delta \left[ g\ast f\left( u\right) \right] ,%
\text{ }t\in \left( 0,T\right) ,\text{ }x\in R^{n},  \tag{1.1}
\end{equation}%
\begin{equation}
u\left( x,0\right) =\varphi \left( x\right) ,\text{ }u_{t}\left( x,0\right)
=\psi \left( x\right) \text{ for a.e. }x\in R^{n},  \tag{1.2}
\end{equation}%
where $A=A\left( x\right) $, $g=g\left( x\right) $ are a linear and
nonlinear operator functions, respectively defined in a Hilbert space $H$; $%
a=a\left( x\right) $ is a complex valued function on $R^{n},$ $f(u)$ is the
given nonlinear function, $\varphi \left( x\right) $ and $\psi \left(
x\right) $ are the given $H-$valued initial functions. The predictions of
classical (local) elasticity theory become inaccurate when the
characteristic length of an elasticity problem is comparable to the atomic
length scale. To solution this situation, a nonlocal theory of elasticity
was introduced (see $[1-3]$ and the references cited therein) and the main
feature of the new theory is the fact that its predictions were more down to
earth than those of the classical theory. For other generalizations of
elasticity we refer the reader to $[4-6]$. \ The global existence of the
Cauchy problem for Boussinesq type nonlocal equations has been studied by
many authors (see $\left[ 7-11\right] $ ). Note that, the existence and
uniqueness of solutions and regularity properties for different type
Boussinesq equations were considered e.g. in $\left[ \text{8-15}\right] $.
Boussinesq type equations occur in a wide variety of physical systems, such
as in the propagation of longitudinal deformation waves in an elastic rod,
hydro-dynamical process in plasma, in materials science which describe
spinodal decomposition and in the absence of mechanical stresses (see $\left[
16-19\right] $). The $L^{p}-$well-posedness of the Cauchy problem $%
(1.1)-\left( 1.2\right) $ depends crucially on the presence of a suitable
kernel. Then the question that naturally arises is which of the possible
forms of the operator functions and kernel functions are relevant for the
global well-posedness of the Cauchy problem $\left( 1.1\right) -\left(
1.2\right) $. In this study, as a partial answer to this question, we
consider the problem $(1.1)-\left( 1.2\right) $ with a general class of
kernel functions with operator coefficients provide local and global
existence for the solutions of $(1.1)-(1.2)$ in frame of $H-$valued $L^{p}$
spaces. The kernel functions most frequently used in the literature are
particular cases of this general class of kernel functions in the scalar
case, i.e. when $H=\mathbb{C}$ ( here, $\mathbb{C}$ denote the set of
complex numbers). In contrast to the above works, we consider the IVP for
nonlocal wave equation with operator coefficients in $H-$valued function
spaces. By selecting the space $H$ and the operators $A$ , $g$ in $\left(
1.1\right) -\left( 1.2\right) ,$ we obtain different classes of nonlocal
wave equations which occur in application. Let we put $H=l_{2}$ and choose $A
$, $g$\ as infinite matrices $\left[ a_{mj}\right] $ and $\left[ g_{mj}%
\right] $, respectively for $m,j=1,2,...,N,$ $N\in \mathbb{N},$ where $%
\mathbb{N-}$denote the set of natural numbers. Then from our results we
obtain the existence, uniqueness and regularity properties of Cauchy problem
for infinity many system of nonlocal WEs 
\begin{equation}
\partial _{t}^{2}u_{m}-a\ast \Delta u_{m}+\dsum\limits_{j=1}^{m}a_{mj}\ast
u_{m}=  \tag{1.3}
\end{equation}%
\[
\dsum\limits_{j=1}^{m}\Delta g_{mj}u_{m}\ast f_{m}\left(
u_{1},u_{2},...,u_{m}\right) ,\text{ }t\in \left[ 0,T\right] \text{, }x\in
R^{n},
\]%
\[
u_{m}\left( x,0\right) =\varphi _{m}\left( x\right) ,\text{ }\partial
_{t}u_{m}\left( x,0\right) =\psi _{m}\left( x\right) ,\text{ }
\]%
where $a_{mj}=a_{mj}\left( x\right) $,$\ g_{mj}=\left( x\right) $ are
complex valued functions, $f_{m}$ are nonlinear functions and $%
u_{j}=u_{j}\left( x,t\right) .$

Moreover, let we choose $E=L^{2}\left( 0,1\right) $ and $A$ to be
degenerated differential operator in $L^{2}\left( 0,1\right) $ defined by 
\[
D\left( A\right) =\left\{ u\in W_{\gamma }^{\left[ 2\right] ,2}\left(
0,1\right) ,\right. \left. \dsum\limits_{i=0}^{\nu _{k}}\alpha _{ki}u^{\left[
i\right] }\left( 0\right) +\beta _{ki}u^{\left[ i\right] }\left( 1\right) =0,%
\text{ }k=1,2\right\} ,\text{ } 
\]%
\begin{equation}
\text{ }A\left( x\right) u=b_{1}\left( x,y\right) u^{\left[ 2\right]
}+b_{2}\left( x,y\right) u^{\left[ 1\right] }\text{, }x\in R^{n},\text{ }%
y\in \left( 0,1\right) ,\text{ }\nu _{k}\in \left\{ 0,1\right\} ,  \tag{1.4}
\end{equation}%
\ \ \ where $u^{\left[ i\right] }=\left( y^{\gamma }\frac{d}{dy}\right)
^{\gamma }u$ for $0\leq \gamma <\frac{1}{2},$ $b_{1}=b_{1}\left( x,y\right) $
is a cont\i nous, $b_{2}=b_{2}\left( x,y\right) $ is a bounded functon on $%
y\in $ $\left[ 0,1\right] $ for a.e. $x\in R^{n},$ $\alpha _{ki}$, $\beta
_{ki}$ are complex numbers and $W_{\gamma }^{\left[ 2\right] ,2}\left(
0,1\right) $ is a weighted Sobolev spase defined by 
\[
W_{\gamma }^{\left[ 2\right] }\left( 0,1\right) =\left\{ {}\right. u:u\in
L^{2}\left( 0,1\right) ,\text{ }u^{\left[ 2\right] }\in L^{2}\left(
0,1\right) ,\text{ } 
\]%
\[
\left\Vert u\right\Vert _{W_{\gamma }^{\left[ 2\right] }}=\left\Vert
u\right\Vert _{L^{2}}+\left\Vert u^{\left[ 2\right] }\right\Vert
_{L^{2}}<\infty . 
\]%
Then, from $\left( 1.1\right) -\left( 1.2\right) $ we get the following
mixed problem for degenerate nonlocal WE

\begin{equation}
u_{tt}-a\ast \Delta u+\left( b_{1}\frac{\partial ^{\left[ 2\right] }u}{%
\partial y^{2}}+b_{2}\frac{\partial ^{\left[ 1\right] }u}{\partial y}\right)
\ast u=\Delta g\ast f\left( u\right) ,\text{ }  \tag{1.5}
\end{equation}%
\[
x\in R^{n},\text{ }y\in \left( 0,1\right) ,\text{ }t\in \left( 0,T\right) ,%
\text{ }u=u\left( x,y,t\right) ,
\]%
\ \ \ 

\begin{equation}
\dsum\limits_{i=0}^{\nu _{k}}\alpha _{ki}u^{\left[ i\right] }\left(
x,0,t\right) +\beta _{ki}u^{\left[ i\right] }\left( x,1,t\right) =0,\text{ }%
k=1,2,  \tag{1.6}
\end{equation}

\begin{equation}
u\left( x,y,0\right) =\varphi \left( x,y\right) ,\text{ }u_{t}\left(
x,y,0\right) =\psi \left( x,y\right) \text{.}  \tag{1.7}
\end{equation}

Note that, the IVP for abstract hyperbolic equations were studied e.g. $%
\left[ \text{20, 21}\right] .$

The strategy is to express the equation $\left( 1.1\right) $ as an integral
equation. To treat the nonlinearity as a small perturbation of the linear
part of the equation, the contraction mapping theorem is used. Also, a
priori estimates on $L^{p}$ norms of solutions of the linearized version are
utilized. The key step is the derivation of the uniform estimate for
solutions of the linearized nonlocal wave equation. The methods of harmonic
analysis, operator theory, interpolation of Banach Spaces and embedding
theorems in Sobolev spaces are the main tools implemented to carry out the
analysis.

In order to state our results precisely, we introduce some notations and
some function spaces.

\begin{center}
\textbf{Definitions and} \textbf{Background}
\end{center}

Let $E$ be a Banach space. $L^{p}\left( \Omega ;E\right) $ denotes the space
of strongly measurable $E$-valued functions that are defined on the
measurable subset $\Omega \subset R^{n}$ with the norm

\[
\left\Vert f\right\Vert _{p}=\left\Vert f\right\Vert _{L^{p}\left( \Omega
;E\right) }=\left( \int\limits_{\Omega }\left\Vert f\left( x\right)
\right\Vert _{E}^{p}dx\right) ^{\frac{1}{p}},\text{ }1\leq p<\infty ,\text{ }
\]

\[
\left\Vert f\right\Vert _{L^{\infty }\left( \Omega ;E\right) }\ =\text{ess}%
\sup\limits_{x\in \Omega }\left\Vert f\left( x\right) \right\Vert _{E}. 
\]

Let $H$ be a Hilbert space. For $p=2$ and $E=H$ the space $L^{p}\left(
\Omega ;E\right) $ becomes the $H$-valued Hilbert space $L^{2}\left( \Omega
;H\right) $ with inner product:%
\[
\left( f,g\right) _{L^{2}\left( \Omega ;H\right) }=\int\limits_{\Omega
}\left( f\left( x\right) ,g\left( x\right) \right) _{H}dx\text{, for any }f,%
\text{ }g\in L^{2}\left( \Omega ;H\right) . 
\]%
For $p=2$ the norm of $L^{p}\left( R^{n};H\right) $ will be denoted just by $%
\left\Vert .\right\Vert _{2}.$

Let $E_{1}$ and $E_{2}$ be two Banach spaces. $\left( E_{1},E_{2}\right)
_{\theta ,p}$ for $\theta \in \left( 0,1\right) ,$ $p\in \left[ 1,\infty %
\right] $ denotes the real interpolation spaces defined by $K$-method $\left[
\text{23, \S 1.3.2}\right] $. Let $E_{1}$ and $E_{2}$ be two Banach spaces. $%
B\left( E_{1},E_{2}\right) $ will denote the space of all bounded linear
operators from $E_{1}$ to $E_{2}.$ For $E_{1}=E_{2}=E$ it will be denoted by 
$B\left( E\right) .$

Here, \ 
\[
S_{\phi }=\left\{ \lambda \in \mathbb{C}\text{, }\left\vert \arg \lambda
\right\vert \leq \phi ,\text{ }0\leq \phi <\pi \right\} .
\]

A closed linear operator\ $A$ is said to be sectorial in a Banach\ space $E$
with bound $M>0$ if $D\left( A\right) $ and $R\left( A\right) $ are dense on 
$E,$ $N\left( A\right) =\left\{ 0\right\} $ and 
\[
\left\Vert \left( A+\lambda I\right) ^{-1}\right\Vert _{B\left( E\right)
}\leq M\left\vert \lambda \right\vert ^{-1} 
\]%
for any $\lambda \in S_{\phi },$ $0\leq \phi <\pi ,$ where $I$ is the
identity operator in $E,$ $B\left( E\right) $ is the space of bounded linear
operators in $E;$ $D\left( A\right) $ and $R\left( A\right) $ denote domain
and range of the operator $A.$ It is known that (see e.g.$\left[ \text{23, 
\S 1.15.1}\right] $) there exist the fractional powers\ $A^{\theta }$ of a
sectorial operator $A.$ Let $E\left( A^{\theta }\right) $ denote the space $%
D\left( A^{\theta }\right) $ with the graphical norm 
\[
\left\Vert u\right\Vert _{E\left( A^{\theta }\right) }=\left( \left\Vert
u\right\Vert ^{p}+\left\Vert A^{\theta }u\right\Vert ^{p}\right) ^{\frac{1}{p%
}},1\leq p<\infty ,\text{ }0<\theta <\infty . 
\]%
A sectorial operator $A\left( \xi \right) $ for $\xi \in R^{n}$ is said to
be uniformly sectorial in a Banach space $E,$ if $D\left( A\left( \xi
\right) \right) $ is independent of $\xi $ and the uniform estimate 
\[
\left\Vert \left( A+\lambda I\right) ^{-1}\right\Vert _{B\left( E\right)
}\leq M\left\vert \lambda \right\vert ^{-1} 
\]%
holds for any $\lambda \in S_{\phi }.$

A uniformly sectorial operator\ $A=A\left( \xi \right) $ belongs to $\sigma
\left( M_{0},\omega ,E\right) $ (see $\left[ \text{29}\right] $ \S\ 11.2) if 
$D\left( A\right) $ is dense on $E,$ $D\left( A\left( \xi \right) \right) $
is independent of $\xi \in R^{n}$ and for $\func{Re}\lambda >\omega $ the
uniform estimate holds 
\[
\left\Vert \left( A\left( \xi \right) -\lambda ^{2}I\right) ^{-1}\right\Vert
_{B\left( E\right) }\leq M_{0}\left\vert \func{Re}\lambda -\omega
\right\vert ^{-1}\text{. }
\]

\textbf{Remark 1.1. }It is known (see e.g. $\left[ \text{30, \S\ 1.6}\right] 
$, Theorem 6.3)\ that if $A\in \sigma \left( M_{0},\omega ,E\right) $ and $%
0\leq \alpha <1$ then it is generates a bounded group operator $U_{A}\left(
t\right) $ satisfying 
\begin{equation}
\left\Vert U_{A}\left( t\right) \right\Vert _{B\left( E\right) }\leq
Me^{\omega \left\vert t\right\vert },\text{ }\left\Vert A^{\alpha
}U_{A}\left( t\right) \right\Vert _{B\left( E\right) }\leq M\left\vert
t\right\vert ^{-\alpha },\text{ }t\in \left( -\infty ,\infty \right) . 
\tag{2.1}
\end{equation}

Let $E$ be a Banach space. $S=S(R^{n};E)$ denotes $E$-valued Schwartz class,
i.e. the space of all $E-$valued rapidly decreasing smooth functions on $%
R^{n}$ equipped with its usual topology generated by seminorms. $S(R^{n};%
\mathbb{C})$ denoted by $S$.

Let $S^{\prime }(R^{n};E)$ denote the space of all continuous linear
operators, $L:S\rightarrow E$, equipped with the bounded convergence
topology. Recall $S(R^{n};E)$ is norm dense in $L^{p}(R^{n};E)$ when $1\leq
p<\infty .$

\ Let $m$ be a positive integer. $W^{m,p}\left( \Omega ;E\right) $ denotes
an $E-$valued Sobolev space, i.e. space of all functions $u\in L^{p}\left(
\Omega ;E\right) $ that have the generalized derivatives $\frac{\partial
^{m}u}{\partial x_{k}^{m}}\in L^{p}\left( \Omega ;E\right) ,$ $1\leq p\leq
\infty $ with the norm 
\[
\ \left\Vert u\right\Vert _{W^{m,p}\left( \Omega ;E\right) }=\left\Vert
u\right\Vert _{L^{p}\left( \Omega ;E\right)
}+\sum\limits_{k=1}^{n}\left\Vert \frac{\partial ^{m}u}{\partial x_{k}^{m}}%
\right\Vert _{L^{p}\left( \Omega ;E\right) }<\infty . 
\]%
\ \ $\ \ $

Let $W^{s,p}\left( R^{n};E\right) $ denotes the fractional Sobolev space of
order for $s\in \mathbb{R}$ that is defined as: 
\[
W^{s,p}\left( E\right) =W^{s,p}\left( R^{n};E\right) =\left\{ u\in \right.
S^{\prime }(R^{n};E), 
\]%
\[
\left\Vert u\right\Vert _{W^{s,p}\left( E\right) }=\left\Vert F^{-1}\left(
I+\left\vert \xi \right\vert ^{2}\right) ^{\frac{s}{2}}\hat{u}\right\Vert
_{L^{p}\left( R^{n};E\right) }<\infty \left. {}\right\} . 
\]%
It clear that $W^{0,p}\left( R^{n};E\right) =L^{p}\left( R^{n};E\right) .$
For $p=2$ and $H$ is a Hilbert space, $W^{s,p}\left( R^{n};H\right) $ will
be denoted just by $H^{s}.$

Let $E_{0}$ and $E$ be two Banach spaces and $E_{0}$ is continuously and
densely embedded into $E$. Here, $W^{s,p}\left( R^{n};E_{0},E\right) $
denote the Sobolev-Lions type space i.e., 
\[
W^{s,p}\left( R^{n};E_{0},E\right) =\left\{ u\in W^{s,p}\left(
R^{n};E\right) \cap L^{p}\left( R^{n};E_{0}\right) ,\right. \text{ } 
\]%
\[
\left. \left\Vert u\right\Vert _{W^{s,p}\left( R^{n};E_{0},E\right)
}=\left\Vert u\right\Vert _{L^{p}\left( R^{n};E_{0}\right) }+\left\Vert
u\right\Vert _{W^{s,p}\left( R^{n};E\right) }<\infty \right\} . 
\]

Let $1\leq p\leq q<\infty .$ A function $\Psi \in L^{\infty }(R^{n})$ is
called a Fourier multiplier from $L^{p}(R^{n};E)$ to $L^{q}(R^{n};E)$ if the
map $P:$ $u\rightarrow F^{-1}\Psi (\xi )Fu$ for $u\in S(R^{n};E)$ is well
defined and extends to a bounded linear operator

\[
P:L^{p}(R^{n};E)\rightarrow L^{q}(R^{n};E). 
\]

Sometimes we use one and the same symbol $C$ without distinction in order to
denote positive constants which may differ from each other even in a single
context. When we want to specify the dependence of such a constant on a
parameter, say $\alpha $, we write $C_{\alpha }$. Moreover, for $u$, $%
\upsilon >0$ the relations $u\lesssim \upsilon ,$ \ $u$\ $\approx $ $%
\upsilon $ means that there exist positive constants $C,$ $C_{1},$ $C_{2}$ \
independent on $u$ and $\upsilon $ such that, respectively 
\[
u\leq C\upsilon ,\text{ }C_{1}\upsilon \leq u\leq C_{2}\upsilon . 
\]

The paper is organized as follows: In Section 1, some definitions and
background are given. In Section 2, we obtain the existence of unique
solution and a priory estimates for solution of the linearized problem $%
(1.1)-\left( 1.2\right) .$ In Section 3, we show the existence and
uniqueness of local strong solution of the problem $(1.1)-\left( 1.2\right) $%
. In the Section 4, we show the same applications of the problem $%
(1.1)-\left( 1.2\right) .$

Sometimes we use one and the same symbol $C$ without distinction in order to
denote positive constants which may differ from each other even in a single
context. When we want to specify the dependence of such a constant on a
parameter, say $h$, we write $C_{h}$.

\begin{center}
\textbf{2. Estimates for linearized equation}
\end{center}

In this section, we make the necessary estimates for solutions of the Cauchy
problem for the nonlocal linear WE 
\begin{equation}
u_{tt}-a\ast \Delta u+A\ast u=g\left( x,t\right) ,\text{ }x\in R^{n},\text{ }%
t\in \left( 0,T\right) ,\text{ }T\in \left( 0,\left. \infty \right] ,\right. 
\tag{2.1}
\end{equation}%
\begin{equation}
u\left( x,0\right) =\varphi \left( x\right) ,\text{ }u_{t}\left( x,0\right)
=\psi \left( x\right) \text{ for a.e. }x\in R^{n},  \tag{2.2}
\end{equation}%
where $A=A\left( x\right) $ is a linear operator function \ defined in a
Hilbert space $H$ and $\ a\geq 0,$

Let $A$ be a sectorial operator in $H.$ Here, 
\[
X_{p}=L^{p}\left( R^{n};H\right) \text{, }X_{p}\left( A^{\gamma }\right)
=L^{p}\left( R^{n};H\left( A^{\gamma }\right) \right) ,\text{ }%
Y_{q}^{s,p}=W^{s,p}\left( R^{n};H\right) \cap X_{q}\text{, } 
\]%
\[
\left\Vert u\right\Vert _{Y_{q}^{s,p}}=\left\Vert u\right\Vert
_{W^{s,p}\left( R^{n};H\right) }+\left\Vert u\right\Vert _{X_{q}}<\infty ,%
\text{ }0<\gamma \leq 1, 
\]%
\[
W^{s,p}\left( A^{\gamma }\right) =W^{s,p}\left( R^{n};H\left( A^{\gamma
}\right) \right) \text{, }Y_{q}^{s,p}\left( A\right) =W^{s,p}\left( A\right)
\cap X_{q}\left( A\right) \text{, } 
\]%
\[
Y^{s,p}\left( A,H\right) =W^{s,p}\left( R^{n};H\left( A\right) ,H\right) ,%
\text{ }Y_{q}^{s,p}\left( A;H\right) =Y^{s,p}\left( A,H\right) \cap X_{q}, 
\]

\[
\left\Vert u\right\Vert _{Y_{q}^{s,p}\left( A;H\right) }=\left\Vert
u\right\Vert _{Y^{s,p}\left( A,H\right) }+\left\Vert u\right\Vert
_{X_{q}}<\infty ,\text{ }1\leq p,\text{ }q\leq \infty . 
\]

\bigskip Let $\hat{A}\left( \xi \right) $\ be the Fourier transformation of $%
A\left( x\right) ,$ i.e. $\hat{A}\left( \xi \right) =F\left( A\left(
x\right) \right) .$ We assume that $\hat{A}\left( \xi \right) $ is uniformly
sectorial operator in a Hilbert\ space $H.$ Let%
\[
\eta =\eta \left( \xi \right) =\left[ \hat{a}\left( \xi \right) \left\vert
\xi \right\vert ^{2}+\hat{A}\left( \xi \right) \right] ^{\frac{1}{2}}.
\]

Let $A$ be a generator of a strongly continuous cosine operator function in
a Hilbert space $H$ defined by formula%
\[
C\left( t\right) =\frac{1}{2}\left( e^{itA^{\frac{1}{2}}}+e^{-itA^{\frac{1}{2%
}}}\right) 
\]%
(see e.g. $\left[ \text{29, \S 11.2, 11.4}\right] $, or $\left[ \text{30}%
\right] $). Then, from the definition of sine operator-function $S\left(
t\right) $ we have%
\[
S\left( t\right) u=\dint\limits_{0}^{t}C\left( \sigma \right) ud\sigma \text{%
, i.e. }S\left( t\right) u=\frac{1}{2i}A^{-\frac{1}{2}}\left( e^{itA^{\frac{1%
}{2}}}-e^{-itA^{\frac{1}{2}}}\right) . 
\]

Let 
\begin{equation}
C\left( t\right) =C\left( \xi ,t\right) =\frac{1}{2}\left( e^{it\eta \left(
\xi \right) }+e^{-it\eta \left( \xi \right) }\right) ,  \tag{2.3}
\end{equation}

\[
\text{ }S\left( t\right) =S\left( \xi ,t\right) =\frac{1}{2i}\eta
^{-1}\left( \xi \right) \left( e^{it\eta \left( \xi \right) }-e^{-it\eta
\left( \xi \right) }\right) .
\]%
\textbf{Condition 2.1. }Assume: (1) $a\in L^{1}\left( R^{n}\right) ,$ $\hat{a%
}\left( \xi \right) \in S\left( \psi \right) $ for all $\xi \in R^{n}$ and $%
\eta \left( \xi \right) \neq 0$\ for all $\xi \in R^{n};$ (2) $\hat{A}\left(
\xi \right) $ is an uniformly sectorial operator in $H$ such that $\hat{A}%
\left( \xi \right) \in \sigma \left( M_{0},\omega ,H\right) ;$ (3) $\varphi
\in W^{s,p}\left( A\right) $ and $\psi \in W^{s,p}\left( A^{\frac{1}{2}%
}\right) ;$ (4)$\ \hat{A}\left( \xi \right) $ is a differentiable operator
function with independent of $\xi $ domain $D\left( D^{\alpha }\hat{A}\left(
\xi \right) \right) =D\left( \hat{A}\right) =$ $D\left( A\right) $ for $%
\alpha =\left( \alpha _{1},\alpha _{2},...,\alpha _{n}\right) ,$ $\left\vert
\alpha \right\vert \leq n$ \ and the uniform estimate holds 
\[
\left\Vert \left[ D^{\alpha }\hat{A}\left( \xi \right) \right] \eta
^{-1}\left( \xi \right) \right\Vert _{B\left( H\right) }\leq M.
\]

First we need the following lemmas:

\textbf{Lemma 2.1. }Let the Assumption (1) of Condition 2.1. holds. Then,
problem $\left( 2.1\right) -\left( 2.2\right) $ has a unique solution.

\textbf{Proof. }By using of the Fourier transform, we get from $(2.1)-\left(
2.2\right) $:%
\begin{equation}
\hat{u}_{tt}\left( \xi ,t\right) +\eta ^{2}\left( \xi \right) \hat{u}\left(
\xi ,t\right) =\hat{g}\left( \xi ,t\right) ,\text{ }  \tag{2.4}
\end{equation}%
\[
\hat{u}\left( \xi ,0\right) =\hat{\varphi}\left( \xi \right) ,\text{ }\hat{u}%
_{t}\left( \xi ,0\right) =\hat{\psi}\left( \xi \right) , 
\]%
where $\hat{u}\left( \xi ,t\right) $ is a Fourier transform of $u\left(
x,t\right) $ with respect to $x$ and $\hat{\varphi}\left( \xi \right) ,$ $%
\hat{\psi}\left( \xi \right) $ are Fourier transform of $\varphi $ and $\psi
,$ respectively. By virtue of $\left[ \text{29, \S 11.2,4}\right] $ we
obtain that $\eta \left( \xi \right) $ is a generator of a strongly
continuous cosine operator function and problem $(2.4)$ has a unique
solution for all $\xi \in R^{n}$ exspressing as%
\begin{equation}
\hat{u}\left( \xi ,t\right) =C\left( \xi ,t\right) \hat{\varphi}\left( \xi
\right) +S\left( \xi ,t\right) \hat{\psi}\left( \xi \right)
+\dint\limits_{0}^{t}S\left( \xi ,t-\tau \right) \hat{g}\left( \xi ,\tau
\right) d\tau ,  \tag{2.5}
\end{equation}%
i.e. problem $(2.1)-\left( 2.2\right) $ has a unique solution 
\begin{equation}
u\left( x,t\right) =C_{1}\left( t\right) \varphi +S_{1}\left( t\right) \psi
+Qg,  \tag{2.6}
\end{equation}%
where $C_{1}\left( t\right) ,$ $S_{1}\left( t\right) ,$ $Q$ are linear
operator functions defined by 
\[
C_{1}\left( t\right) \varphi =F^{-1}\left[ C\left( \xi ,t\right) \hat{\varphi%
}\left( \xi \right) \right] ,\text{ }S_{1}\left( t\right) \psi =F^{-1}\left[
S\left( \xi ,t\right) \hat{\psi}\left( \xi \right) \right] , 
\]

\[
Qg=F^{-1}\tilde{Q}\left( \xi ,t\right) ,\text{ }\tilde{Q}\left( \xi
,t\right) =\dint\limits_{0}^{t}F^{-1}\left[ S\left( \xi ,t-\tau \right) \hat{%
g}\left( \xi ,\tau \right) \right] d\tau . 
\]%
\textbf{Theorem 2.1. }Let the Condition 2.1 holds and $s>1+\frac{n}{p}$ with 
$p\in \left( 1,\infty \right) $. Then for $\varphi \in Y_{1}^{s,p}\left(
A;H\right) $ and $\psi \in Y_{1}^{s,p}\left( A^{\frac{1}{2}};H\right) ,$ $%
g\left( x,t\right) \in Y_{1}^{s,p}$ problem $(2.1)-(2.2)$ has a unique
generalized solution 
\[
u(x,t)\in C^{2}\left( \left[ 0,T\right] ;Y^{s}\left( A;H\right) \right) . 
\]%
Moreover, the following estimate holds 
\begin{equation}
\left\Vert A^{\frac{1}{2}}u\right\Vert _{X_{\infty }}+\left\Vert A^{\frac{1}{%
2}}u_{t}\right\Vert _{X_{\infty }}\leq C\left[ \left\Vert \varphi
\right\Vert _{Y_{1}^{s,p}\left( A\right) }\right. +  \tag{2.7}
\end{equation}

\[
\left\Vert \psi \right\Vert _{Y_{1}^{s,p}\left( A^{\frac{1}{2}}\right)
}+\left. \dint\limits_{0}^{t}\left( \left\Vert g\left( .,\tau \right)
\right\Vert _{Y^{s,p}}+\left\Vert g\left( .,\tau \right) \right\Vert
_{X_{1}}\right) d\tau \right] , 
\]%
uniformly in $t\in \left[ 0,T\right] $, where the positive constant $C$
depends only on initial data and the space $H$.

\textbf{Proof. }From Lemma 2.1 we obtain that the problem $(2.1)-(2.2)$ has
a unique generalized solution $u(x,t)\in C^{2}\left( \left[ 0,T\right]
;Y^{s,p}\left( A;H\right) \right) $ for $\varphi ,$ $\psi \in
Y_{1}^{s,p}\left( A\right) ,$ $g\left( x,t\right) \in Y_{1}^{s,p}$. Let $%
N\in \mathbb{N}$ and 
\[
\Pi _{N}=\left\{ \xi :\xi \in R^{n},\text{ }\left\vert \xi \right\vert \leq
N\right\} ,\text{ }\Pi _{N}^{\prime }=\left\{ \xi :\xi \in R^{n},\text{ }%
\left\vert \xi \right\vert \geq N\right\} . 
\]

From $\left( 2.6\right) $ we deduced that

\[
\left\Vert A^{\frac{1}{2}}u\right\Vert _{X_{\infty }}\lesssim \left\Vert
F^{-1}C\left( \xi ,t\right) \hat{A}^{\frac{1}{2}}\hat{\varphi}\left( \xi
\right) \right\Vert _{L^{\infty }\left( \Pi _{N}\right) }+ 
\]%
\begin{equation}
\left\Vert F^{-1}S\left( \xi ,t\right) \hat{A}^{\frac{1}{2}}\hat{\psi}\left(
\xi \right) \right\Vert _{L^{\infty }\left( \Pi _{N}\right) }+\left\Vert
F^{-1}C\left( \xi ,t\right) \hat{A}^{\frac{1}{2}}\hat{\varphi}\left( \xi
\right) \right\Vert _{L^{\infty }\left( \Pi _{N}^{\prime }\right) }+ 
\tag{2.8}
\end{equation}%
\[
\left\Vert F^{-1}S\left( \xi ,t\right) \hat{A}^{\frac{1}{2}}\hat{\psi}\left(
\xi \right) \right\Vert _{L^{\infty }\left( \Pi _{N}^{\prime }\right) }+%
\frac{1}{2}\left\Vert F^{-1}\hat{A}^{\frac{1}{2}}\tilde{Q}\left( \xi
,t\right) \right\Vert _{L^{\infty }\left( \Pi _{N}\right) }+ 
\]%
\[
\frac{1}{2}\left\Vert F^{-1}\hat{A}^{\frac{1}{2}}\tilde{Q}\left( \xi
,t\right) \hat{g}\left( \xi ,\tau \right) \right\Vert _{L^{\infty }\left(
\Pi _{N}^{\prime }\right) }. 
\]

\bigskip Due to uniform boundedness of operator functions $C\left( \xi
,t\right) $, $S\left( \xi ,t\right) $, in view of $\left( 2.3\right) $ and\
by Minkowski's inequality for integrals\ we get 
\[
\left\Vert F^{-1}C\left( \xi ,t\right) \hat{A}^{\frac{1}{2}}\hat{\varphi}%
\left( \xi \right) \right\Vert _{L^{\infty }\left( \Pi _{N}\right)
}+\left\Vert \left\Vert F^{-1}S\left( \xi ,t\right) \hat{A}^{\frac{1}{2}}%
\hat{\psi}\left( \xi \right) \right\Vert \right\Vert _{L^{\infty }\left( \Pi
_{N}\right) }\lesssim 
\]

\begin{equation}
\left[ \left\Vert A\varphi \right\Vert _{X_{1}}+\left\Vert A^{\frac{1}{2}%
}\psi \right\Vert _{X}+\left\Vert g\right\Vert _{X_{1}}\right] .  \tag{2.9}
\end{equation}%
Moreover, from $\left( 2.6\right) $ we deduced that 
\[
\left\Vert F^{-1}C\left( \xi ,t\right) \hat{A}^{\frac{1}{2}}\hat{\varphi}%
\left( \xi \right) \right\Vert _{L^{\infty }\left( \Pi _{N}^{\prime }\right)
}+\left\Vert F^{-1}S\left( \xi ,t\right) \hat{A}^{\frac{1}{2}}\hat{\psi}%
\left( \xi \right) \right\Vert _{L^{\infty }\left( \Pi _{N}^{\prime }\right)
}\lesssim 
\]%
\[
\left\Vert F^{-1}C\left( \xi ,t\right) \hat{A}^{\frac{1}{2}}\hat{\varphi}%
\left( \xi \right) \right\Vert _{L^{\infty }\left( \Pi _{N}^{\prime }\right)
}+\left\Vert F^{-1}S\left( \xi ,t\right) \hat{A}^{\frac{1}{2}}\hat{\psi}%
\left( \xi \right) \right\Vert _{L^{\infty }\left( \Pi _{N}^{\prime }\right)
}+ 
\]%
\[
\left\Vert F^{-1}S\left( \xi ,t\right) \hat{A}^{\frac{1}{2}}\tilde{Q}\left(
\xi ,t\right) \hat{g}\left( \xi ,\tau \right) \right\Vert _{L^{\infty
}\left( \Pi _{N}^{\prime }\right) }\lesssim 
\]%
\begin{equation}
\left\Vert F^{-1}\left( 1+\left\vert \xi \right\vert ^{2}\right) ^{-\frac{s}{%
2}}C\left( \xi ,t\right) \left( 1+\left\vert \xi \right\vert ^{2}\right) ^{%
\frac{s}{2}}\hat{A}^{\frac{1}{2}}\hat{\varphi}\left( \xi \right) \right\Vert
_{L^{\infty }\left( \Pi _{N}^{\prime }\right) }+  \tag{2.10}
\end{equation}%
\[
\left\Vert F^{-1}\left( 1+\left\vert \xi \right\vert ^{2}\right) ^{-\frac{s}{%
2}}S\left( \xi ,t\right) \left( 1+\left\vert \xi \right\vert ^{2}\right) ^{%
\frac{s}{2}}\hat{A}^{\frac{1}{2}}\hat{\psi}\left( \xi \right) \right\Vert
_{L^{\infty }\left( \Pi _{N}^{\prime }\right) }+ 
\]%
\[
\left\Vert F^{-1}\left( 1+\left\vert \xi \right\vert ^{2}\right) ^{-\frac{s}{%
2}}S\left( \xi ,t\right) \left( 1+\left\vert \xi \right\vert ^{2}\right) ^{%
\frac{s}{2}}\hat{A}^{\frac{1}{2}}\tilde{Q}\left( \xi ,t\right) \hat{g}\left(
\xi ,\tau \right) \right\Vert _{L^{\infty }\left( \Pi _{N}^{\prime }\right)
}; 
\]%
here, the space $L^{\infty }\left( \Omega ;H\right) $ was denoted by $%
L^{\infty }\left( \Omega \right) $. From $\left( 2.3\right) $ it clear to
see that%
\[
\frac{\partial }{\partial \xi _{k}}\left[ \left( 1+\left\vert \xi
\right\vert ^{2}\right) ^{-\frac{s}{2}}C\left( \xi ,t\right) \right] =-s\xi
_{k}\left( 1+\left\vert \xi \right\vert ^{2}\right) ^{-\frac{s}{2}-1}C\left(
\xi ,t\right) + 
\]%
\begin{equation}
\frac{t}{4}\left( 1+\left\vert \xi \right\vert ^{2}\right) ^{-\frac{s}{2}%
}\eta ^{\frac{1}{2}}\left( \xi \right) \left( 2\xi _{k}a+\frac{\partial }{%
\partial \xi _{k}}\hat{A}\left( \xi \right) \right) S\left( \xi ,t\right) . 
\tag{2.11}
\end{equation}%
By assumption (4) and in view of $s>1+\frac{n}{p}$\ from $\left( 2.3\right) $%
, $\left( 2.11\right) $\ we obtain%
\[
\sup\limits_{\xi \in R^{n},t\in \left[ 0,T\right] }\left\vert \xi
\right\vert ^{\left\vert \alpha \right\vert +\frac{n}{p}}\left\Vert
D^{\alpha }\left[ \left( 1+\left\vert \xi \right\vert ^{2}\right) ^{-\frac{s%
}{2}}C\left( \xi ,t\right) \right] \right\Vert _{B\left( H\right) }\leq
C_{1}, 
\]%
\ 
\begin{equation}
\sup\limits_{\xi \in R^{n},t\in \left[ 0,T\right] }\left\vert \xi
\right\vert ^{\left\vert \alpha \right\vert +\frac{n}{p}}\left\Vert
D^{\alpha }\left[ \left( 1+\left\vert \xi \right\vert ^{2}\right) ^{-\frac{s%
}{2}}S\left( \xi ,t\right) \right] \right\Vert _{B\left( H\right) }\leq C_{2}
\tag{2.12}
\end{equation}%
for $\alpha =\left( \alpha _{1},\alpha _{2},...,\alpha _{n}\right) $, $%
\alpha _{k}\in \left\{ 0,1\right\} $ and $\xi \in R^{n}$ uniformly in $t\in %
\left[ 0,T\right] .$ Hence, by Fourier multiplier theorems (see e.g. $\left[ 
\text{22, Theorem 4.3}\right] $), from $\left( 2.12\right) $ we get that the
functions $\left( 1+\left\vert \xi \right\vert ^{2}\right) ^{-\frac{s}{2}%
}C\left( \xi ,t\right) ,$ $\left( 1+\left\vert \xi \right\vert ^{2}\right)
^{-\frac{s}{2}}S\left( \xi ,t\right) $ are $L^{p}\left( R^{n};H\right)
\rightarrow L^{\infty }\left( R^{n};H\right) $ Fourier multipliers. Then by
Minkowski's inequality for integrals, from $\left( 2.3\right) ,$ $\left(
2.10\right) $ and $\left( 2.11\right) -\left( 2.12\right) $ we have%
\[
\left\Vert F^{-1}C\left( \xi ,t\right) \hat{A}^{\frac{1}{2}}\hat{\varphi}%
\left( \xi \right) \right\Vert _{L^{\infty }\left( \Pi _{N}^{\prime }\right)
}+\left\Vert F^{-1}S\left( \xi ,t\right) \hat{A}^{\frac{1}{2}}\hat{\psi}%
\left( \xi \right) \right\Vert _{L^{\infty }\left( \Pi _{N}^{\prime }\right)
}\lesssim 
\]

\begin{equation}
\left[ \left\Vert A\varphi \right\Vert _{H^{s,p}}+\left\Vert A^{\frac{1}{2}%
}\psi \right\Vert _{H^{s,p}}+\left\Vert g\right\Vert _{H^{s,p}}\right] . 
\tag{2.13}
\end{equation}%
By reasoning as the above, we have 
\begin{equation}
\left\Vert F^{-1}\hat{A}^{\frac{1}{2}}\tilde{Q}\left( \xi ,t\right)
\right\Vert _{X_{\infty }}\leq C\dint\limits_{0}^{t}\left( \left\Vert
g\left( .,\tau \right) \right\Vert _{H^{s,p}}+\left\Vert g\left( .,\tau
\right) \right\Vert _{X_{1}}\right) d\tau  \tag{2.14}
\end{equation}%
uniformly in $t\in \left[ 0,T\right] $. Thus, from $\left( 2.6\right) ,$ $%
\left( 2.13\right) $ and $\left( 2.14\right) $ we obtain 
\begin{equation}
\left\Vert A^{\frac{1}{2}}u\right\Vert _{X_{\infty }}\leq C\left[ \left\Vert
A^{\frac{1}{2}}\varphi \right\Vert _{Y^{s,p}}+\left\Vert A^{\frac{1}{2}%
}\varphi \right\Vert _{X_{1}}\right. +  \tag{2.15}
\end{equation}

\[
\left\Vert A^{\frac{1}{2}}\psi \right\Vert _{Y^{s,p}}+\left\Vert A^{\frac{1}{%
2}}\psi \right\Vert _{X_{1}}+\left. \dint\limits_{0}^{t}\left( \left\Vert
g\left( .,\tau \right) \right\Vert _{Y^{s,p}}+\left\Vert g\left( .,\tau
\right) \right\Vert _{X_{1}}\right) d\tau \right] . 
\]%
By differentiating $\left( 2.6\right) $, in a similar way we obtain 
\begin{equation}
\left\Vert A^{\frac{1}{2}}u_{t}\right\Vert _{X_{\infty }}\leq C\left[
\left\Vert A\varphi \right\Vert _{Y^{s,p}}+\left\Vert A\varphi \right\Vert
_{X_{1}}\right. +  \tag{2.16}
\end{equation}

\[
\left\Vert A^{\frac{1}{2}}\psi \right\Vert _{Y^{s,p}}+\left\Vert A^{\frac{1}{%
2}}\psi \right\Vert _{X_{1}}+\left. \dint\limits_{0}^{t}\left( \left\Vert
g\left( .,\tau \right) \right\Vert _{Y^{s,p}}+\left\Vert g\left( .,\tau
\right) \right\Vert _{X_{1}}\right) d\tau \right] . 
\]

Then from $\left( 2.15\right) $\ and $\left( 2.16\right) $ we obtain the
assertion.

\textbf{Theorem 2.2. }Let the Condition 2.1 holds and $s>1+\frac{n}{p}$.
Then for $g\left( x,t\right) \in W^{s,p}$ the solution of $\left( 2.1\right)
-\left( 2.2\right) $ satisfies the following uniform estimate%
\begin{equation}
\left( \left\Vert A^{\frac{1}{2}}u\right\Vert _{H^{s,p}}+\left\Vert A^{\frac{%
1}{2}}u_{t}\right\Vert _{H^{s,p}}\right) \leq  \tag{2.17}
\end{equation}

\[
C_{0}\left( \left\Vert A\varphi \right\Vert _{H^{s,p}}+\left\Vert A^{\frac{1%
}{2}}\psi \right\Vert _{H^{s,p}}+\dint\limits_{0}^{t}\left\Vert g\left(
.,\tau \right) \right\Vert _{H^{s,p}}d\tau \right) . 
\]%
\textbf{Proof. }From $\left( 2.5\right) $ and $\left( 2.11\right) $ we get
the following uniform estimate 
\begin{equation}
\left( \left\Vert F^{-1}\left( 1+\left\vert \xi \right\vert ^{2}\right) ^{%
\frac{s}{2}}\hat{A}^{\frac{1}{2}}\hat{u}\right\Vert _{X_{p}}+\left\Vert
F^{-1}\left( 1+\left\vert \xi \right\vert ^{2}\right) ^{\frac{s}{2}}\hat{A}^{%
\frac{1}{2}}\hat{u}_{t}\right\Vert _{X_{p}}\right) \leq  \tag{2.18}
\end{equation}

\[
C\left\{ \left\Vert F^{-1}\left( 1+\left\vert \xi \right\vert ^{2}\right) ^{%
\frac{s}{2}}C\left( \xi ,t\right) \hat{A}^{\frac{1}{2}}\hat{\varphi}%
\right\Vert _{X_{p}}\right. +\left\Vert F^{-1}\left( 1+\left\vert \xi
\right\vert ^{2}\right) ^{\frac{s}{2}}\hat{A}^{\frac{1}{2}}S\left( \xi
,t\right) \hat{\psi}\right\Vert _{X_{p}}+ 
\]

\[
\left. \dint\limits_{0}^{t}\left\Vert \left( 1+\left\vert \xi \right\vert
^{2}\right) ^{\frac{s}{2}}\hat{A}^{\frac{1}{2}}\tilde{Q}\left( \xi ,t\right) 
\hat{g}\left( \xi ,\tau \right) \right\Vert _{X_{p}}d\tau \right\} . 
\]

\bigskip By Condition 2.1 and by usingFourier multiplier theorem $\left[ 
\text{22, Theorem 4.3}\right] $ and by reasoning as in Theorem 2.1\ we get
that $C\left( \xi ,t\right) $, $S\left( \xi ,t\right) $ and $\hat{A}^{\frac{1%
}{2}}S\left( \xi ,t\right) $\ are Fourier multipliers in $L^{p}\left(
R^{n};H\right) $ uniformly with respect to $t\in \left[ 0,T\right] .$ So,
the estimate $\left( 2.18\right) $ by using the Minkowski's inequality for
integrals implies $\left( 2.17\right) .$

\begin{center}
\textbf{3. Local well posedness of IVP for nonlinear nonlocal WE}
\end{center}

In this section, we will show the local existence and uniqueness of solution
for the Cauchy problem $(1.1)-(1.2).$ For the study of the nonlinear problem 
$\left( 1.1\right) -\left( 1.2\right) $ we need the following lemmas

\textbf{Lemma 3.1} (Abstract Nirenberg's inequality). Let $H$ be a Hilbert
space. Assume that $u\in L_{p}\left( \Omega ;H\right) $, $D^{m}u$ $\in
L_{q}\left( \Omega ;H\right) $, $p,q\in \left( 1,\infty \right) $. Then for $%
i$ with $0\leq i\leq m,$ $m>\frac{n}{q}$ we have 
\begin{equation}
\left\Vert D^{i}u\right\Vert _{r}\leq C\left\Vert u\right\Vert _{p}^{1-\mu
}\dsum\limits_{k=1}^{n}\left\Vert D_{k}^{m}u\right\Vert _{q}^{\mu }, 
\tag{3.1}
\end{equation}%
where%
\[
\frac{1}{r}=\frac{i}{m}+\mu \left( \frac{1}{q}-\frac{m}{n}\right) +\left(
1-\mu \right) \frac{1}{p},\text{ }\frac{i}{m}\leq \mu \leq 1. 
\]

\textbf{Proof. }By virtue of interpolation of Banach spaces $\left[ \text{%
23, \S 1.3.2}\right] ,$ in order to prove $\left( 3.1\right) $ for any given 
$i,$ one has only to prove it for the extreme values $\mu =\frac{i}{m}$ and $%
\mu =1$. For the case of $\mu =1$, i.e., $\frac{1}{r}=\frac{i}{m}+\frac{1}{q}%
-\frac{m}{n}$ the estimate $\left( 3.1\right) $ is obtained from Theorem A$%
_{1}$. The case $\mu =\frac{i}{m}$ is derived by reasoning as in $\left[ 
\text{25, \S\ 2 }\right] $ and in replacing absolute value of complex-valued
function $u$ by the $H-$norm of $H$-valued function.

\bigskip Note that, for $H=\mathbb{C}$ the lemma considered by L. Nirenberg $%
\left[ 25\right] .$

Using the chain rule of the composite function, from Lemma 3.1 we can prove
the following result

\textbf{Lemma 3.2. }Let $H$ be a Hilbert space. Assume that $u\in $ $%
W^{m,p}\left( \Omega ;H\right) \cap L^{\infty }\left( \Omega ;H\right) $,
and $f\left( u\right) $ possesses continuous derivatives up to order $m\geq
1 $. Then $f\left( u\right) -f\left( 0\right) \in W^{m,p}\left( \Omega
;H\right) $ and 
\[
\left\Vert f\left( u\right) -f\left( 0\right) \right\Vert _{p}\leq
\left\Vert f^{^{\left( 1\right) }}\left( u\right) \right\Vert _{\infty
}\left\Vert u\right\Vert _{p}, 
\]

\begin{equation}
\left\Vert D^{k}f\left( u\right) \right\Vert _{p}\leq
C_{0}\dsum\limits_{j=1}^{k}\left\Vert f^{\left( j\right) }\left( u\right)
\right\Vert _{\infty }\left\Vert u\right\Vert _{\infty }^{j-1}\left\Vert
D^{k}u\right\Vert _{p}\text{, }1\leq k\leq m,  \tag{3.2}
\end{equation}%
where $C_{0}$ $\geq 1$ is a constant.

For $H=\mathbb{C}$ the lemma coincide with the corresponding inequality in $%
\left[ 26\right] .$ Let $H_{0}$ denotes the real interpolation space between 
$Y^{s,p}\left( A,H\right) $ and $L^{p}\left( R^{n};H\right) $ with $\theta =%
\frac{1}{2p}$, i.e. 
\[
\text{ }H_{0}=\left( Y^{s,p}\left( A,H\right) ,L^{p}\left( R^{n};H\right)
\right) _{\frac{1}{2p},p}. 
\]%
\textbf{Remark 3.1. }By using J.Lions-I. Petree result (see e.g $\left[ 
\text{21, \S\ 1.8.}\right] $) we obtain that the map $u\rightarrow u\left(
t_{0}\right) $, $t_{0}\in \left[ 0,T\right] $ is continuous and surjective
from $H^{s}$ onto $H_{0}$ and there is a constant $C_{1}$ such that 
\[
\left\Vert u\left( t_{0}\right) \right\Vert _{H_{0}}\leq C_{1}\left\Vert
u\right\Vert _{W^{s,p}\left( R^{n};H\right) },\text{ }1\leq p\leq \infty 
\text{.} 
\]%
First all of, we define the space $Y\left( T\right) =C\left( \left[ 0,T%
\right] ;Y_{\infty }^{s,p}\left( A,H\right) \right) $ equipped with the norm
defined by%
\[
\left\Vert u\right\Vert _{Y\left( T\right) }=\max\limits_{t\in \left[ 0,T%
\right] }\left\Vert u\right\Vert _{Y_{\infty }^{s,p}\left( A,H\right) },%
\text{ }u\in Y\left( T\right) . 
\]

It is easy to see that $Y\left( T\right) $ is a Banach space. For $\varphi $%
, $\psi \in Y_{\infty }^{s,p}\left( A^{\frac{1}{2}}\right) $, let 
\[
M=\left\Vert \varphi \right\Vert _{Y_{\infty }^{s,p}\left( A^{\frac{1}{2}%
}\right) }+\left\Vert \psi \right\Vert _{Y_{\infty }^{s,p}\left( A^{\frac{1}{%
2}}\right) }. 
\]

\textbf{Definition 3.1. }For any $T>0$, $\varphi ,$ $\psi \in Y_{\infty
}^{s,p}\left( A^{\frac{1}{2}}\right) $, the function $u$ $\in C^{2}\left( %
\left[ 0,T\right] ;Y_{\infty }^{s,p}\left( A,H\right) \right) $ satisfies
the equation $(1.1)-(1.2)$ is called the continuous solution\ or the strong
solution of the problem $(1.1)-(1.2).$ If $T<\infty $, then $u\left(
x,t\right) $ is called the local strong solution of the problem $%
(1.1)-(1.2). $ If $T=\infty $, then $u\left( x,t\right) $ is called the
global strong solution of $(1.1)-(1.2)$.

\textbf{Condition 3.1. }Assume:

(1) the Condition 2.1 holds for $s>\frac{n}{p}$ and $\varphi \in Y_{\infty
}^{s,p}\left( A\right) ,$ $\psi \in Y_{\infty }^{s,p}\left( A^{\frac{1}{2}%
}\right) $;

(2) the kernel $g=g\left( x\right) $ is a bounded integrable operator
function in $H,$ whose Fourier transform satisfies 
\[
0\leq \left\Vert \hat{g}\left( \xi \right) \right\Vert _{B\left( H\right)
}\lesssim \left( 1+\left\vert \xi \right\vert ^{2}\right) ^{-1}\text{ for
all }\xi \in R^{n}; 
\]

(3) the function $u\rightarrow $ $f\left( x,t,u\right) $: $R^{n}\times \left[
0,T\right] \times H_{0}\rightarrow H$ is a measurable in $\left( x,t\right)
\in R^{n}\times \left[ 0,T\right] $ for $u\in H_{0}.$ Moreover, $f\left(
x,t,u\right) $ is continuous in $u\in H_{0}$ and $f\left( x,t,u\right) \in
C^{\left[ s\right] +1}\left( H_{0};H\right) $ uniformly with respect to $%
x\in R^{n},$ $t\in \left[ 0,T\right] .$

Main aim of this section is to prove the following result:

\textbf{Theorem 3.1. }Let the Condition 3.1 holds. Then problem $\left(
1.1\right) -\left( 1.2\right) $ has a unique local strange solution $u\in
C^{\left( 2\right) }\left( \left[ 0,\right. \left. T_{0}\right) ;Y_{\infty
}^{s,p}\left( A,H\right) \right) $, where $T_{0}$ is a maximal time interval
that is appropriately small relative to $M$. Moreover, if

\begin{equation}
\sup_{t\in \left[ 0\right. ,\left. T_{0}\right) }\left( \left\Vert
u\right\Vert _{Y_{\infty }^{s,p}\left( A;H\right) }+\left\Vert
u_{t}\right\Vert _{Y_{\infty }^{s,p}\left( A;H\right) }\right) <\infty 
\tag{3.3}
\end{equation}%
then $T_{0}=\infty .$

\textbf{Proof. }First, we are going to prove the existence and the
uniqueness of the local continuous solution of $(1.1)-\left( 1.2\right) $ by
contraction mapping principle. By $(2.5)$, $\left( \left( 2.6\right) \right) 
$ the problem of finding a solution $u$ of $(1.1)-\left( 1.2\right) $ is
equivalent to finding a fixed point of the mapping

\[
G\left( u\right) =C_{1}\left( t\right) \varphi \left( x\right) +S_{1}\left(
t\right) \psi \left( x\right) +Q\left( u\right) ,
\]%
where $C_{1}\left( t\right) ,$ $S_{1}\left( t\right) $ are defined by $%
\left( 2.6\right) $ and $Q\left( u\right) $ is a map in $Y\left( T,M\right) $
defined by 
\begin{equation}
Q\left( u\right) =-i\dint\limits_{0}^{t}F^{-1}\left[ U\left( \xi ,t-\tau
\right) \left\vert \xi \right\vert ^{2}\hat{g}\left( \xi \right) \hat{f}%
\left( u\right) \left( \xi ,\tau \right) \right] d\tau ,  \tag{3.4}
\end{equation}%
where%
\[
Y\left( T;M\right) =\left\{ u:u\in L^{q}\left( \left[ 0,T\right]
;L^{r}\left( R^{n};H\left( A\right) \right) \right) ,\left\Vert
Au\right\Vert _{_{L_{t}^{q}L_{x}^{r}\left( H\right) }}\leq M\right\} \text{ }
\]%
with $T$ and $M$ to be determined. So, we will find $T$ and $M$ so that $G$
is a contraction on $Y(T,M)$.

From Lemma 3.2 we know that $f(u)\in $ $L^{p}\left( 0,T;Y_{\infty
}^{s,p}\right) $ for any $T>0$. Thus, by Lemma 2.1, problem $(1.1)-\left(
1.2\right) $ has a solution satisfies the following 
\begin{equation}
G\left( u\right) \left( x,t\right) =C_{1}\left( t\right) \varphi
+S_{1}\left( t\right) \psi +Qu,  \tag{3.5}
\end{equation}%
where $C_{1}\left( t\right) $, $S_{1}\left( t\right) $ are defined by $%
\left( 2.5\right) $ and$\left( 2.6\right) .$ From Lemma 3.2 it is easy to
see that the map $G$ is well defined for $f\in C^{\left[ s\right] +1}\left(
H_{0};\mathbb{C}\right) $. First,\ let us prove that the map $G$ has a
unique fixed point in $Q\left( M;T\right) .$ For this aim, it is sufficient
to show that the operator $G$ maps $Q\left( M;T\right) $ into $Q\left(
M;T\right) $ and $G:$ $Q\left( M;T\right) $ $\rightarrow $ $Q\left(
M;T\right) $ is strictly contractive if $T$ is appropriately small relative
to $M.$ Consider the function \ $\bar{f}\left( \sigma \right) $: $\left[
0,\right. $ $\left. \infty \right) \rightarrow \left[ 0,\right. $ $\left.
\infty \right) $ defined by 
\[
\ \bar{f}\left( \sigma \right) =\max\limits_{\left\vert x\right\vert \leq
\sigma }\left\{ \left\Vert f^{\left( 1\right) }\left( x\right) \right\Vert _{%
\mathbb{C}},\left\Vert f^{\left( 2\right) }\left( x\right) \right\Vert _{%
\mathbb{C}}\text{ ,...,}\left\Vert f^{\left[ s\right] }\left( x\right)
\right\Vert _{\mathbb{C}}\right\} ,\text{ }\sigma \geq 0. 
\]

It is clear to see that the function $\bar{f}\left( \sigma \right) $ is
continuous and nondecreasing on $\left[ 0,\right. $ $\left. \infty \right) .$
From Lemma 3.2 we have\qquad

\[
\left\Vert f\left( u\right) \right\Vert _{Y^{s,2}}\leq \left\Vert f^{\left(
1\right) }\left( u\right) \right\Vert _{X_{\infty }}\left\Vert u\right\Vert
+\left\Vert f^{\left( 1\right) }\left( u\right) \right\Vert _{X_{\infty
}}\left\Vert Du\right\Vert + 
\]

\begin{equation}
C_{0}\left[ \left\Vert f^{\left( 1\right) }\left( u\right) \right\Vert
_{X_{\infty }}\left\Vert u\right\Vert +\left\Vert f^{\left( 2\right) }\left(
u\right) \right\Vert _{X_{\infty }}\left\Vert u\right\Vert _{X_{\infty
}}\left\Vert D^{2}u\right\Vert \right] +...+  \tag{3.6}
\end{equation}

\[
\left\Vert f^{\left( \left[ s\right] \right) }\left( u\right) \right\Vert
_{X_{\infty }}\left\Vert u\right\Vert _{X_{\infty }}\left\Vert D^{^{\left[ s%
\right] }}u\right\Vert \leq 2C_{0}\bar{f}\left( M+1\right) \left( M+1\right)
\left\Vert u\right\Vert _{Y^{s,2}}. 
\]%
\ In view of the assumpt\i on (1) and by using Minkowski's inequality for
integrals\"{o} we obtain from $\left( 3.5\right) $:%
\begin{equation}
\left\Vert G\left( u\right) \right\Vert _{X_{\infty }}\lesssim \left\Vert
\varphi \right\Vert _{\infty }+\left\Vert \psi \right\Vert _{\infty
}+\dint\limits_{0}^{t}\left\Vert \Delta \left[ g\ast f\left( \left( u\right)
\right) \right] \left( x,\tau \right) d\tau \right\Vert _{\infty }, 
\tag{3.7}
\end{equation}%
\begin{equation}
\left\Vert G\left( u\right) \right\Vert _{Y^{2,p}}\lesssim \left\Vert
\varphi \right\Vert _{Y^{s,p}}+\left\Vert \psi \right\Vert
_{Y^{s,p}}+\dint\limits_{0}^{t}\left\Vert \Delta \left[ g\ast f\left(
u\right) \right] \left( x,\tau \right) d\tau \right\Vert _{Y^{2,p}}d\tau . 
\tag{3.8}
\end{equation}%
Thus, from $\left( 3.6\right) -\left( 3.8\right) $ and Lemma 3.2 we get 
\[
\left\Vert G\left( u\right) \right\Vert _{Y\left( T\right) }\leq M+T\left(
M+1\right) \left[ 1+2C_{0}\left( M+1\right) \bar{f}\left( M+1\right) \right]
. 
\]%
If $T$ satisfies 
\begin{equation}
T\leq \left\{ \left( M+1\right) \left[ 1+2C_{0}\left( M+1\right) \bar{f}%
\left( M+1\right) \right] \right\} ^{-1},  \tag{3.9}
\end{equation}%
then 
\[
\left\Vert Gu\right\Vert _{Y\left( T\right) }\leq M+1. 
\]%
Therefore, if $\left( 3.9\right) $ holds, then $G$ maps $Q\left( M;T\right) $
into $Q\left( M;T\right) .$ Now, we are going to prove that the map $G$ is
strictly contractive. Assume $T>0$ and $u_{1},$ $u_{2}\in $ $Q\left(
M;T\right) $ given. We get%
\[
G\left( u_{1}\right) -G\left( u_{2}\right) = 
\]%
\[
\dint\limits_{0}^{t}F^{-1}\left[ S\left( t-\tau ,\xi \right) \left\vert \xi
\right\vert ^{2}\hat{g}\left( \xi \right) \left( \hat{f}\left( u_{1}\right)
\left( \xi ,\tau \right) -\hat{f}\left( u_{2}\right) \left( \xi ,\tau
\right) \right) \right] d\tau ,\text{ }t\in \left( 0,T\right) . 
\]%
By using the assumption (3) and the mean value theorem, we obtain%
\[
\hat{f}\left( u_{1}\right) -\hat{f}\left( u_{2}\right) =\hat{f}^{\left(
1\right) }\left( u_{2}+\eta _{1}\left( u_{1}-u_{2}\right) \right) \left(
u_{1}-u_{2}\right) ,\text{ } 
\]

\[
D_{\xi }\left[ \hat{f}\left( u_{1}\right) -\hat{f}\left( u_{2}\right) \right]
=\hat{f}^{\left( 2\right) }\left( u_{2}+\eta _{2}\left( u_{1}-u_{2}\right)
\right) \left( u_{1}-u_{2}\right) D_{\xi }u_{1}+\text{ } 
\]%
\[
\hat{f}^{\left( 1\right) }\left( u_{2}\right) \left( D_{\xi }u_{1}-D_{\xi
}u_{2}\right) , 
\]%
\[
D_{\xi }^{2}\left[ \hat{f}\left( u_{1}\right) -\hat{f}\left( u_{2}\right) %
\right] =\hat{f}^{\left( 3\right) }\left( u_{2}+\eta _{3}\left(
u_{1}-u_{2}\right) \right) \left( u_{1}-u_{2}\right) \left( D_{\xi
}u_{1}\right) ^{2}+\text{ } 
\]%
\[
\hat{f}^{\left( 2\right) }\left( u_{2}\right) \left( D_{\xi }u_{1}-D_{\xi
}u_{2}\right) \left( D_{\xi }u_{1}+D_{\xi }u_{2}\right) + 
\]%
\[
\hat{f}^{\left( 2\right) }\left( u_{2}+\eta _{4}\left( u_{1}-u_{2}\right)
\right) \left( u_{1}-u_{2}\right) D_{\xi }^{2}u_{1}+\hat{f}^{\left( 1\right)
}\left( u_{2}\right) \left( D_{\xi }^{2}u_{1}-D_{\xi }^{2}u_{2}\right) , 
\]%
where $0<\eta _{i}<1,$ $i=1,2,3,4.$ Thus, using Hollder's and Nirenberg's
inequality, we have%
\begin{equation}
\left\Vert \hat{f}\left( u_{1}\right) -\hat{f}\left( u_{2}\right)
\right\Vert _{X_{\infty }}\leq \bar{f}\left( M+1\right) \left\Vert
u_{1}-u_{2}\right\Vert _{X_{\infty }},  \tag{3.10}
\end{equation}%
\begin{equation}
\left\Vert \hat{f}\left( u_{1}\right) -\hat{f}\left( u_{2}\right)
\right\Vert \leq \bar{f}\left( M+1\right) \left\Vert u_{1}-u_{2}\right\Vert ,
\tag{3.11}
\end{equation}%
\begin{equation}
\left\Vert D_{\xi }\left[ \hat{f}\left( u_{1}\right) -\hat{f}\left(
u_{2}\right) \right] \right\Vert \leq \left( M+1\right) \bar{f}\left(
M+1\right) \left\Vert u_{1}-u_{2}\right\Vert _{X_{\infty }}+  \tag{3.12}
\end{equation}%
\[
\bar{f}\left( M+1\right) \left\Vert \hat{f}\left( u_{1}\right) -\hat{f}%
\left( u_{2}\right) \right\Vert ,...,+ 
\]%
\[
\left\Vert D_{\xi }^{\left[ s\right] }\left[ \hat{f}\left( u_{1}\right) -%
\hat{f}\left( u_{2}\right) \right] \right\Vert \leq \left( M+1\right) \bar{f}%
\left( M+1\right) \left\Vert u_{1}-u_{2}\right\Vert _{X_{\infty }}\left\Vert
D_{\xi }^{\left[ s\right] }u_{1}\right\Vert ^{2}+ 
\]%
\[
\bar{f}\left( M+1\right) \left\Vert D_{\xi }\left( u_{1}-u_{2}\right)
\right\Vert _{4}\left\Vert D_{\xi }\left( u_{1}+u_{2}\right) \right\Vert
_{4}+ 
\]%
\[
\bar{f}\left( M+1\right) \left\Vert u_{1}-u_{2}\right\Vert _{X_{\infty
}}\left\Vert D_{\xi }^{2}u_{1}\right\Vert +\bar{f}\left( M+1\right)
\left\Vert D_{\xi }\left( u_{1}-u_{2}\right) \right\Vert \leq 
\]%
\begin{equation}
C^{2}\bar{f}\left( M+1\right) \left\Vert u_{1}-u_{2}\right\Vert _{X_{\infty
}}\left\Vert u_{1}\right\Vert _{X_{\infty }}\left\Vert D_{\xi
}^{2}u_{1}\right\Vert +  \tag{3.13}
\end{equation}%
\[
C^{2}\bar{f}\left( M+1\right) \left\Vert u_{1}-u_{2}\right\Vert _{X_{\infty
}}^{\frac{1}{2}}\left\Vert D_{\xi }^{2}\left( u_{1}-u_{2}\right) \right\Vert
\left\Vert u_{1}+u_{2}\right\Vert _{X_{\infty }}^{\frac{1}{2}}\left\Vert
D_{\xi }^{2}\left( u_{1}+u_{2}\right) \right\Vert + 
\]%
\[
\left( M+1\right) \bar{f}\left( M+1\right) \left\Vert u_{1}-u_{2}\right\Vert
_{X_{\infty }}+\bar{f}\left( M+1\right) \left\Vert D_{\xi }^{2}\left(
u_{1}-u_{2}\right) \right\Vert \leq 
\]%
\[
3C^{2}\left( M+1\right) ^{2}\bar{f}\left( M+1\right) \left\Vert
u_{1}-u_{2}\right\Vert _{X_{\infty }}+2C^{2}\left( M+1\right) \bar{f}\left(
M+1\right) \left\Vert D_{\xi }^{2}\left( u_{1}-u_{2}\right) \right\Vert , 
\]%
where $C$ is the constant in Lemma $3.1$. From $\left( 3.10\right) -\left(
3.11\right) $, using Minkowski's inequality for integrals and Young's
inequality, we obtain%
\[
\left\Vert G\left( u_{1}\right) -G\left( u_{2}\right) \right\Vert _{Y\left(
T\right) }\leq \dint\limits_{0}^{t}\left\Vert u_{1}-u_{2}\right\Vert
_{X_{\infty }}d\tau +\dint\limits_{0}^{t}\left\Vert u_{1}-u_{2}\right\Vert
_{Y^{s,p}}d\tau + 
\]%
\[
\dint\limits_{0}^{t}\left\Vert f\left( u_{1}\right) -f\left( u_{2}\right)
\right\Vert _{X_{\infty }}d\tau +\dint\limits_{0}^{t}\left\Vert f\left(
u_{1}\right) -f\left( u_{2}\right) \right\Vert _{Y^{s,2}}d\tau \leq 
\]%
\[
T\left[ 1+C_{1}\left( M+1\right) ^{2}\bar{f}\left( M+1\right) \right]
\left\Vert u_{1}-u_{2}\right\Vert _{Y\left( T\right) }, 
\]%
where $C_{1}$ is a constant. If $T$ satisfies $\left( 3.9\right) $ and the
following inequality holds 
\begin{equation}
T\leq \frac{1}{2}\left[ 1+C_{1}\left( M+1\right) ^{2}\bar{f}\left(
M+1\right) \right] ^{-1},  \tag{3.14}
\end{equation}%
then 
\[
\left\Vert Gu_{1}-Gu_{2}\right\Vert _{Y\left( T\right) }\leq \frac{1}{2}%
\left\Vert u_{1}-u_{2}\right\Vert _{Y\left( T\right) }. 
\]%
That is, $G$ is a contractive map. By contraction mapping principle we know
that $G(u)$ has a fixed point $u(x,t)\in $ $Q\left( M;T\right) $ that is a
solution of $(1.1)-(1.2)$. From $\left( 2.9\right) -\left( 2.11\right) $ we
get that $u$ is a solution of the following integral equation 
\[
u\left( x,t\right) =C_{1}\left( t\right) \varphi +S_{1}\left( t\right) \psi
+ 
\]%
\[
\dint\limits_{0}^{t}F^{-1}\left[ S\left( t-\tau ,\xi \right) \left\vert \xi
\right\vert ^{2}\hat{g}\left( \xi \right) \hat{f}\left( u\right) \left( \xi
,\tau \right) \right] d\tau ,\text{ }t\in \left( 0,T\right) . 
\]%
Let us show that this solution is a unique in $Y\left( T\right) $. Let $%
u_{1} $, $u_{2}\in Y\left( T\right) $ are two solution of the problem $%
(1.1)-(1.2)$. Then%
\begin{equation}
u_{1}-u_{2}=\dint\limits_{0}^{t}F^{-1}\left[ S\left( t-\tau ,\xi \right)
\left\vert \xi \right\vert ^{2}\hat{g}\left( \xi \right) \left( \hat{f}%
\left( u_{1}\right) \left( \xi ,\tau \right) -\hat{f}\left( u_{2}\right)
\left( \xi ,\tau \right) \right) \right] d\tau .  \tag{3.15}
\end{equation}%
By the definition of the space $Y\left( T\right) $, we can assume that%
\[
\left\Vert u_{1}\right\Vert _{X_{\infty }}\leq C_{1}\left( T\right) ,\text{ }%
\left\Vert u_{1}\right\Vert _{X_{\infty }}\leq C_{1}\left( T\right) . 
\]%
Hence, by Minkowski's inequality for integrals and by Theorem 2.2 \ from $%
\left( 3.15\right) $ we obtain

\begin{equation}
\left\Vert u_{1}-u_{2}\right\Vert _{Y^{s,p}}\leq C_{2}\left( T\right) \text{ 
}\dint\limits_{0}^{t}\left\Vert u_{1}-u_{2}\right\Vert _{Y^{s,p}}d\tau . 
\tag{3.16}
\end{equation}%
From $(3.16)$ and Gronwall's inequality, we have $\left\Vert
u_{1}-u_{2}\right\Vert _{Y^{s,p}}=0$, i.e. problem $(1.1)-(1.2)$ has a
unique solution which belongs to $Y\left( T\right) .$ That is, we obtain the
first part of the assertion. Now, let $\left[ 0\right. ,\left. T_{0}\right) $
be the maximal time interval of existence for $u\in Y\left( T_{0}\right) $.
It remains only to show that if $(3.3)$ is satisfied, then $T_{0}=\infty $.
Assume contrary that, $\left( 3.3\right) $ holds and $T_{0}<\infty .$ For $%
T\in \left[ 0\right. ,\left. T_{0}\right) ,$ we consider the following
integral equation

\begin{equation}
\upsilon \left( x,t\right) =C_{1}\left( t\right) u\left( x,T\right)
+S_{1}\left( t\right) u_{t}\left( x,T\right) -  \tag{3.17}
\end{equation}

\[
\dint\limits_{0}^{t}F^{-1}\left[ S\left( t-\tau ,\xi \right) \left\vert \xi
\right\vert ^{2}\hat{g}\left( \xi \right) \hat{f}\left( \upsilon \right)
\left( \xi ,\tau \right) \right] d\tau ,\text{ }t\in \left( 0,T\right) . 
\]%
By virtue of $(3.3)$, for $T^{\prime }>T$ we have 
\[
\sup_{t\in \left[ 0\right. ,\left. T\right) }\left( \left\Vert u\right\Vert
_{H^{s,p}\left( A\right) }+\left\Vert u\right\Vert _{\infty }+\left\Vert
u_{t}\right\Vert _{H^{s,p}\left( A\right) }+\left\Vert u_{t}\right\Vert
_{\infty }\right) <\infty . 
\]

By reasoning as a first part of theorem and by contraction mapping
principle, there is a $T^{\ast }\in \left( 0,T_{0}\right) $ such that for
each $T\in \left[ 0\right. ,\left. T_{0}\right) ,$ the equation $\left(
3.17\right) $ has a unique solution $\upsilon \in Y\left( T^{\ast }\right) .$
The estimates $\left( 3.9\right) $ and $\left( 3.14\right) $ imply that $%
T^{\ast }$ can be selected independently of $T\in \left[ 0\right. ,\left.
T_{0}\right) .$ Set $T=T_{0}-\frac{T^{\ast }}{2}$ and define 
\begin{equation}
\tilde{u}\left( x,t\right) =\left\{ 
\begin{array}{c}
u\left( x,t\right) ,\text{ }t\in \left[ 0,T\right] \\ 
\upsilon \left( x,t-T\right) \text{, }t\in \left[ T,T_{0}+\frac{T^{\ast }}{2}%
\right]%
\end{array}%
\right. .  \tag{3.18}
\end{equation}%
By construction $\tilde{u}\left( x,t\right) $ is a solution of the problem $%
(1.1)-(1.2)$ on $\left[ T,T_{0}+\frac{T^{\ast }}{2}\right] $ and in view of
local uniqueness, $\tilde{u}\left( x,t\right) $ extends $u.$ This is against
to the maximality of $\left[ 0\right. ,\left. T_{0}\right) $, i.e we obtain $%
T_{0}=\infty .$

Here, we will denote $L^{2}\left( R^{n};H\right) $ by $L^{2}$. Let 
\[
W^{s,p}\left( R^{n};E\right) ,\text{ }W^{s,p}\left( R^{n};E\left( A^{\theta
}\right) \right) 
\]
will be denoted by $H^{s}$, $H^{s}\left( A^{\theta }\right) $ respectively,
for $E=H$ and $p=2.$ First, \ we show the following lemmas concerning the
behaviour of the nonlinear term in $H-$valued space $H^{s},$ in a similar
way as $[$8, 13, 27$]$.

\textbf{\ Lemma 3.3.} Let $s\geq 0,$ $f\in C^{\left[ s\right] +1}\left(
R;H\right) $ with $f(0)=0$. Then for any $u\in H^{s}\cap L^{\infty }$, we
have $f(u)\in H^{s}\cap L^{\infty }.$ Moreover, there is some constant $A(M)$
depending on $M$ such that for all $u\in H^{s}\cap L^{\infty }$ with $%
\left\Vert u\right\Vert _{L^{\infty }}\leq M,$%
\[
\left\Vert f(u)\right\Vert _{H^{s}}\leq A\left( M\right) \left\Vert
u)\right\Vert _{H^{s}}. 
\]%
\textbf{Lemma 3.4. } Let $s\geq 0,$ $f\in C^{\left[ s\right] +1}\left(
R;H\right) $. Then for for any $M$ there is some constant $B(M)$ depending
on $M$ such that for all $u$, $\upsilon \in H^{s}\cap L^{\infty }$ with $%
\left\Vert u\right\Vert _{L^{m}}\leq M,$ $\left\Vert \upsilon \right\Vert
_{L^{\infty }}\leq M,$ $\left\Vert u\right\Vert _{H^{8}}\leq M,$ $\left\Vert
\upsilon \right\Vert _{H^{s}}\leq M,$%
\[
\left\Vert f(u)-f(\upsilon \right\Vert _{H^{s}}\leq B\left( M\right)
\left\Vert u-\upsilon \right\Vert _{H^{s}},\text{ }\left\Vert
f(u)-f(\upsilon \right\Vert _{L^{\infty }}\leq B\left( M\right) \left\Vert
u-\upsilon \right\Vert _{L^{\infty }}. 
\]

By reasoning as in $\left[ \text{13, Lemma 3.4}\right] $ and $\left[ \text{%
28, Lemma X 4}\right] $ we have, respectively

\textbf{Corollary 3.1.} Let $s>\frac{n}{2},$ $f\in C^{\left[ s\right]
+1}\left( R;H\right) $. Then for any $B$ there is a constant $B(M)$
depending on $M$ such that for all $u$, $\upsilon \in H^{s}$ with $%
\left\Vert u\right\Vert _{H^{s}}\leq M,$ $\left\Vert \upsilon \right\Vert
_{H^{s}}\leq M,$%
\[
\left\Vert f(u)-f(\upsilon \right\Vert _{H^{s}}\leq B\left( M\right)
\left\Vert u-\upsilon \right\Vert _{H^{s}}. 
\]

\bigskip \textbf{Lemma 3.5. } If $s>0$, then $Y_{\infty }^{s,2}$ is an
algebra. Moreover, for \ $f,$ $g\in Y_{\infty }^{s,2},$ 
\[
\left\Vert fg\right\Vert _{H^{s}}\leq C\left[ \left\Vert f\right\Vert
_{\infty }+\left\Vert g\right\Vert _{H^{s}}+\left\Vert f\right\Vert
_{H^{s}}+\left\Vert g\right\Vert _{\infty }\right] . 
\]

\textbf{\ }By using Lemmas 3.3, 3.5 we obta\i n

\textbf{Lemma 3.6 .} Let $s\geq 0,$ $f\in C^{\left[ s\right] +1}\left(
R;H\right) $ and $f\left( u\right) =O\left( \left\vert u\right\vert ^{\alpha
+1}\right) $ \ for $u\rightarrow 0$, $\alpha \geq 1$ be a positive integer.
If $u\in Y_{\infty }^{s,2}$ and $\left\Vert u\right\Vert _{\infty }\leq M$,
\ then 
\[
\left\Vert f(u)\right\Vert _{H^{s}}\leq C\left( M\right) \left[ \left\Vert
u\right\Vert _{H^{s}}\left\Vert u\right\Vert _{\infty }^{\alpha }\right] , 
\]%
\[
\left\Vert f(u)\right\Vert _{1}\leq C\left( M\right) \left\Vert u\right\Vert
^{2}\left\Vert u\right\Vert _{\infty }^{\alpha -1}. 
\]

\textbf{Lemma 3.7 }$\left[ \text{13, Lemma 3.4}\right] $\textbf{.} Let $%
s\geq 0,$ $f\in C^{\left[ s\right] +1}\left( R;H\right) $ and $f\left(
u\right) =O\left( \left\vert u\right\vert ^{\alpha +1}\right) $ \ for $%
u\rightarrow 0$, $\alpha \geq 0$ be a positive integer. If $u,$ $\upsilon
\in Y_{\infty }^{s,2},$ $\left\Vert u\right\Vert _{H^{s}}\leq M$, \ $%
\left\Vert \upsilon \right\Vert _{H^{s}}\leq M$ and $\left\Vert u\right\Vert
_{\infty }\leq M$, \ $\left\Vert \upsilon \right\Vert _{\infty }\leq M,$
then 
\[
\left\Vert f(u)-f(\upsilon )\right\Vert _{H^{s}}\leq C\left( M\right) \left[
\left( \left\Vert u\right\Vert _{\infty }-\left\Vert \upsilon \right\Vert
_{\infty }\right) \left( \left\Vert u\right\Vert _{H^{s}}+\left\Vert
\upsilon \right\Vert _{H^{s}}\right) \right. 
\]%
\[
\left( \left\Vert u\right\Vert _{\infty }+\left\Vert \upsilon \right\Vert
_{\infty }\right) ^{\alpha -1}, 
\]%
\[
\left\Vert f(u)-f(\upsilon \right\Vert _{1}\leq C\left( M\right) \left(
\left\Vert u\right\Vert _{\infty }+\left\Vert \upsilon \right\Vert _{\infty
}\right) ^{\alpha -1}\left( \left\Vert u\right\Vert _{2}+\left\Vert \upsilon
\right\Vert _{2}\right) \left\Vert u-\upsilon \right\Vert _{2}. 
\]

Consider the problem $\left( 1.1\right) -\left( 1.2\right) ,$ when $\varphi
, $ $\psi \in H^{s}.$ $\ $By reasoning as in Theorem 3.1 and $\left[ \text{%
13, Theorem 1.1}\right] $ we have:

\textbf{Condition 3.2. }Assume: (1) the Condition 2.1 holds; (2) $\varphi
\in H^{s}\left( A\right) ,$ $\psi \in H^{s}\left( A^{\frac{1}{2}}\right) $ $%
\ $and $s>\frac{n}{2};$ (3) $f\in C^{\left[ s\right] }\left( R;H\right) $
with $f(0)=0$; (4) the kernel $g$ is a bounded integrable operator function
in $H,$ whose Fourier transform satisfies 
\begin{equation}
0\leq \left\Vert \hat{g}\left( \xi \right) \right\Vert _{B\left( H\right)
}\leq C_{g}\left( 1+\left\vert \xi \right\vert ^{2}\right) ^{-\frac{r}{2}}%
\text{ for all }\xi \in R^{n}\text{ and }r\geq 2.  \tag{3.19}
\end{equation}

\ \textbf{Theorem 3.2. }Let the Condition 3.2 holds. Assume $f\in
C^{k}\left( R;H\right) ,$\ with $k$ an integer $k\geq s>\frac{n}{2},$
satisfies $f\left( u\right) =O\left( \left\vert u\right\vert ^{\alpha
+1}\right) $ \ for $u\rightarrow 0.$ Then there exists a constant $\delta
>0, $ such that for any $\varphi ,$ $\psi \in Y_{1}^{s,2}\left( A^{\frac{1}{2%
}}\right) $ satisfying 
\begin{equation}
\left\Vert \varphi \right\Vert _{Y_{1}^{s,2}\left( A^{\frac{1}{2}}\right)
}+\left\Vert \psi \right\Vert _{Y_{1}^{s,2}\left( A^{\frac{1}{2}}\right)
}\leq \delta ,  \tag{3.20}
\end{equation}%
problem $\left( 1.1\right) -\left( 1.2\right) $ has a unique local strange
solution $u\in C^{\left( 2\right) }\left( \left[ 0,\right. \left. \infty
\right) ;Y^{s,2}\left( A,H\right) \right) $. Moreover,

\begin{equation}
\sup_{0\leq t<\infty }\left( \left\Vert u\right\Vert _{Y^{s,2}\left( A^{%
\frac{1}{2}}\right) }+\left\Vert u_{t}\right\Vert _{Y^{s,2}\left( A^{\frac{1%
}{2}}\right) }\right) \leq C\delta ,  \tag{3.21}
\end{equation}%
where the constant $C$ only depends on $f$ \ and initial data.

\textbf{Proof. }Consider a metric space \ defined by 
\[
W^{s}=\left\{ u\in ^{\left( 2\right) }\left( \left[ 0,\right. \left. \infty
\right) ;Y^{s,2}\left( A,H\right) \right) ,\text{ }\left\Vert u\right\Vert
_{W^{s}}\leq 3C_{0}\delta \right\} , 
\]%
equipped with the norm%
\[
\left\Vert u\right\Vert _{W^{s}}=\sup\limits_{t\geq 0}\left( \left\Vert
u\right\Vert _{Y_{\infty }^{s,2}\left( A;H\right) }+\left\Vert
u_{t}\right\Vert _{Y_{\infty }^{s,2}\left( A;H\right) }\right) , 
\]%
where $\delta >0$ \ satisfies $\left( 3.20\right) $ and $C_{0}$ \ is a
constant in Theorem 2.1. \ It is easy to prove that $W^{s}$ is a complete
metric space. From Sobolev imbedding theorem we know that $\left\Vert
u\right\Vert _{\infty }\leq 1$ if we take that $\delta $ is enough small.
Consider the problem $\left( 3.4\right) $. From Lemma 3.6 we get that $%
f\left( u\right) \in L^{2}\left( 0,T;Y_{1}^{s,2}\right) $ for any $T>0$.
Thus the problem $\left( 3.4\right) $ has a unique solution which can be
written as $\left( 3.5\right) .$ We should prove that the operator $G\left(
u\right) $\ defined by $\left( 3.5\right) $ is strictly contractive if $%
\delta $ is suitable small. \ In fact, by $(2.17)$ in Theorem 2.2 and Lemma
3.6 we get 
\[
\left\Vert A^{\frac{1}{2}}G\left( u\right) \right\Vert _{H^{s}}+\left\Vert
A^{\frac{1}{2}}G_{t}\left( u\right) \right\Vert _{H^{s}}\leq C_{0}\left[
\left\Vert A^{\frac{1}{2}}\varphi \right\Vert _{H^{s}}\right. +\left\Vert A^{%
\frac{1}{2}}\psi \right\Vert _{H^{s}}+ 
\]

\[
\left. \dint\limits_{0}^{t}\left\Vert K\left( u\right) \left( .,\tau \right)
\right\Vert _{H^{s}}d\tau \right] \leq C_{0}\delta
+C_{0}\dint\limits_{0}^{t}\left\Vert K\left( u\right) \left( .,\tau \right)
\right\Vert _{H^{s}}d\tau \leq 
\]%
\begin{equation}
C_{0}\delta +C\dint\limits_{0}^{t}\left\Vert u\left( \tau \right)
\right\Vert _{H^{s}}\left\Vert u\left( \tau \right) \right\Vert _{\infty
}^{\alpha }d\tau \leq C_{0}\delta +C\left\Vert u\right\Vert _{W^{s}}^{\alpha
+1},  \tag{3.22}
\end{equation}%
where 
\[
K\left( u\right) \left( .,\tau \right) =F^{-1}\left[ S\left( t-\tau ,\xi
\right) \left\vert \xi \right\vert ^{2}\hat{g}\left( \xi \right) \hat{f}%
\left( u\right) \left( \xi ,\tau \right) \right] . 
\]%
Therefore, from $(3.22)$ we have 
\begin{equation}
\left\Vert G\left( u\right) \right\Vert _{W^{s}}\leq 2C_{0}\delta
+C\left\Vert u\right\Vert _{W^{s}}^{\alpha +1}\text{.}  \tag{3.23}
\end{equation}%
Taking that $\delta $ is enough small such that $C\left( 3C_{0}\delta
\right) ^{\alpha }$ $<1/3$, from $\left( 3.23\right) $ and from Theorems
2.1, 2.2 we get that $G$ maps $W^{s}$ into $W^{s}$. Then, by reasoning as
the remaining part of $\left[ \text{13, Theorem 1.1}\right] $ we obtain that 
$G$ :$W^{s}\rightarrow W^{s}$ is strictly contractive. Using the contraction
mapping principle, we know that $G(u)$ has a unique fixed point $u(x,t)\in $ 
$C^{\left( 2\right) }([0,\infty );Y^{s,2}\left( A,H\right) )$ and $u(x,t)$
is the solution of the problem $(1.1)-(1.2)$. Moreover, by virtue of Theorem
2.1 from $\left( 3.20\right) $ we obtain $\left( 3.21\right) .$

We claim that the solution $u(x,t)$ of the problem $(1.1)-(1.2)$ is also
unique in $C^{1}([0,\infty );Y^{s,2}\left( A,H\right) )$. In fact, let $u_{1}
$ and $u_{2}$ be two solutions of the problem $(1.1)-(1.2)$ and $u_{1}$, $%
u_{2}\in C^{\left( 2\right) }([0,\infty );Y^{s,2}\left( A,H\right) )$. Let $%
u=u_{1}-u_{2}$; then%
\[
u_{tt}-a\ast \Delta u+A\ast u=\Delta \left[ g\ast \left( f\left(
u_{1}\right) -f\left( u_{2}\right) \right) \right] .
\]

\bigskip This fact is derived in a similar way as in Theorem 3.2, by using
Theorems 2.1, 2.2 and Gronwall's inequality

\textbf{Theorem 3.3. }Let the Condition 3.2 hold. Then there is some $T>0$
such that the problem $(1.1)-(1.2)$ for initial data $\varphi ,$ $\psi \in
H^{s}$ is well posed with solution in $C^{1}\left( \left[ 0,T\right]
;Y^{s,2}\left( A,H\right) \right) .$

\textbf{Proof. }Consider the convolution operator $u\rightarrow \Delta \left[
g\ast f\left( u\right) \right] .$ In view of assumptions we have \ 
\[
\left\Vert \Delta g\ast \upsilon \right\Vert _{H^{s}}\lesssim \left\Vert
\left( 1+\xi \right) ^{\frac{s}{2}}\left\vert \xi \right\vert ^{2}\hat{g}%
\left( \xi \right) \hat{\upsilon}\left( \xi \right) \right\Vert \lesssim
\left\Vert \upsilon \right\Vert _{H^{s}},
\]%
i.e. $\Delta g\ast \upsilon $ is a bounded linear operator on $H^{s}.$ Then
by Corollary 3.1, $K\left( u\right) $ is locally Lipschitz on $H^{s}$. Then
by reasoning as in Theorem 3.2 and $\left[ \text{13, Theorem 1.1}\right] $
we obtain that $G$: $H^{s}\rightarrow H^{s}$ is strictly contractive. Using
the contraction mapping principle, we get that the operator $G(u)$ defined
by $\left( 3.5\right) $ has a unique fixed point $u(x,t)\in $ $C^{\left(
2\right) }([0,\infty );Y^{s,2}\left( A,H\right) )$ \ and $u(x,t)$ is the
solution of the problem $(1.1)-(1.2)$. Moreover, we show that the solution $%
u(x,t)$ of $(1.1)-(1.2)$ is also unique in $C^{\left( 2\right) }([0,\infty
);Y^{s,2}\left( A,H\right) )$. In fact, let $u_{1}$ and $u_{2}$ be two
solutions of the problem $(1.1)-(1.2)$ and $u_{1}$, $u_{2}\in C^{\left(
2\right) }([0,\infty );Y^{s,2}\left( A,H\right) )$. Let $u=u_{1}-u_{2}$; then%
\[
u_{tt}-a\ast \Delta u+A\ast u=\Delta \left[ g\ast \left( f\left(
u_{1}\right) -f\left( u_{2}\right) \right) \right] .
\]%
This fact is derived in a similar way as in Theorem 3.2, by using Theorems
2.1, 2.2 and Gronwall's inequality.

\textbf{Theorem 3.4. }Let the Condition 3.2 holds for $r>2+\frac{n}{2}$.
Then there is some $T>0$ such that the problem $(1.1)-(1.2)$ is well posed
for $\varphi ,$ $\psi \in Y_{\infty }^{s,2}$ with solution in $C^{\left(
2\right) }\left( \left[ 0,T\right] ;Y_{\infty }^{s,2}\left( A;H\right)
\right) .$

\textbf{Proof. } All we need here, is to show that $K\ast f(u)$ is Lipschitz
on $Y_{\infty }^{s,2}$. Indeed, by reasoning as in Theorem 3.3 we have 
\[
\left\Vert \Delta g\ast \upsilon \right\Vert _{H^{s+r-2}}\lesssim \left\Vert
\left( 1+\left\vert \xi \right\vert ^{2}\right) ^{\frac{s+r-2}{2}}\left\vert
\xi \right\vert ^{2}\hat{g}\left( \xi \right) \hat{\upsilon}\left( \xi
\right) \right\Vert \lesssim \left\Vert \upsilon \right\Vert _{H^{s}}, 
\]

\bigskip Then $\Delta g\ast \upsilon $ is a bounded linear map from $H^{s}$
into $H^{s+r-2}$. Since $s\geq 0$ and $r$ $>2+\frac{n}{2}$ \ we get $s+r-2>%
\frac{n}{2}.$ The embedding theorem for $H-$valued Sobolev spaces (see e.g, $%
\left[ 31\right] $) implies that $\Delta g\ast \upsilon $\ is a bounded
linear map from $Y_{\infty }^{s,2}\left( A;H\right) $ into $Y_{\infty
}^{s,2}\left( A;H\right) $. Lemma 3.4 implies the Lipschitz condition on $%
Y_{\infty }^{s,2}$. Then, by reasoning as in Theorem 3.3 we obtain the
assertion.

The solution in theorems 3.2-3.4 can be extended to a maximal interval $%
[0,T_{\max }),$ where finite $T_{\max }$ is characterized by the blow-up
condition 
\[
\limsup\limits_{T\rightarrow T_{\max }}\left\Vert u\right\Vert _{Y_{\infty
}^{s,2}\left( A;H\right) }=\infty . 
\]

\textbf{Lemma 3.8.} \textbf{\ }Let the Condition 3.2 hold and $u$ is a
solution of $(1.1)-(1.2).$ Then there is a global solution if for any $%
T<\infty $ we have%
\begin{equation}
\sup_{t\in \left[ 0,\right. \left. T\right) }\left( \left\Vert u\right\Vert
_{Y_{\infty }^{s,2}\left( A;H\right) }+\left\Vert u_{t}\right\Vert
_{Y_{\infty }^{s,2}\left( A;H\right) }\right) <\infty .  \tag{3.24}
\end{equation}

\textbf{Proof. }Indeed, by reasoning as in the second part of the proof of
Theorem 3.1, by using a continuation of local solution of $(1.1)-(1.2)$ and
assuming contrary that, $\left( 3.24\right) $ holds and $T_{0}<\infty $\ \
we obtain contradiction, i.e. we get $T_{0}=T_{\max }=\infty .$

\begin{center}
\textbf{4. Conservation of energy and global existence. }
\end{center}

In this section, we prove the existence and the uniqueness of the global
strong solution for the problem $(1.1)-(1.2).$ \ For this purpose, we are
going to make a priori estimates of the local strong solution of $%
(1.1)-(1.2).$

\textbf{Condition 4.1. }Suppose the Condition 3.2 is satisfied.\textbf{\ }%
Moreover\textbf{, }assume that the kernel $g$ is a bounded operator function
in $H,$ whose Fourier transform satisfies%
\[
0<\left\Vert \hat{g}\left( \xi \right) \right\Vert _{B\left( H\right)
}\lesssim \left( 1+\left\vert \xi \right\vert ^{2}\right) ^{-\frac{r}{2}}%
\text{ for all }\xi \in R^{n}\text{ and }r\leq 2\left( s+1\right) . 
\]

Let $F^{-1}$ denote the inverse Fourier transform. We consider the operator $%
B$ defined by 
\[
u\in D\left( B\right) =H^{s},\text{ }Bu=F^{-1}\left[ \left\vert \xi
\right\vert ^{-1}\hat{g}^{^{-\frac{1}{2}}}\left( \xi \right) \hat{u}\left(
\xi \right) \right] , 
\]%
Then it is clear to see that 
\begin{equation}
B^{-2}u=-\Delta g\ast u,\text{ }B^{-1}u=F^{-1}\left[ \left\vert \xi
\right\vert \hat{g}^{\frac{1}{2}}\left( \xi \right) \hat{u}\left( \xi
\right) \right] .  \tag{4.1}
\end{equation}

First, we show the following\bigskip

\textbf{Lemma 4.1. }Let the Condition 4.1 holds. Assume the solution of $%
(1.1)-(1.2)$ exists in $C^{\left( 2\right) }\left( \left[ 0,T\right]
;Y_{\infty }^{s,2}\left( A;H\right) \right) $. Then%
\[
\hat{A}^{\frac{1}{2}}Bu,\text{ }\hat{A}^{\frac{1}{2}}Bu_{t}\in C^{\left(
1\right) }\left( \left[ 0,\right. \left. T\right) ;L^{2}\right) . 
\]

\bigskip \textbf{Proof.} By Lemma 2.1, problem $\left( 1.1\right) -\left(
1.2\right) $ is equ\i valent to following integral equation ,%
\begin{equation}
u\left( x,t\right) =C_{1}\left( t\right) \varphi +S_{1}\left( t\right) \psi
+Qg,  \tag{4.2}
\end{equation}%
where $C_{1}\left( t\right) $, $S_{1}\left( t\right) $ are are operator
functions defined by $\left( 2.5\right) $ and $\left( 2.6\right) $, $\ $%
where $g$ replaced by $g\ast f\left( u\right) $ and 
\begin{equation}
Qg=\dint\limits_{0}^{t}F^{-1}\left[ S\left( \xi ,t-\tau \right) \left\vert
\xi \right\vert ^{2}\hat{g}\left( \xi \right) \hat{f}\left( u\right) \left(
\xi \right) \right] d\tau .  \tag{4.3}
\end{equation}%
From $\left( 4.2\right) $ we get that 
\[
u_{t}\left( x,t\right) =\frac{d}{dt}C_{1}\left( t\right) \varphi +\frac{d}{dt%
}S_{1}\left( t\right) \psi + 
\]%
\begin{equation}
\dint\limits_{0}^{t}F^{-1}\left[ C\left( \xi ,t-\tau \right) \left\vert \xi
\right\vert ^{2}\hat{g}\left( \xi \right) \hat{f}\left( G\left( u\right)
\left( \xi \right) \right) \right] d\tau .  \tag{4.4}
\end{equation}%
Since $C_{1}\left( t\right) ,$ $S_{1}\left( t\right) $ and $\frac{d}{dt}%
S\left( \xi ,t\right) $\ are uniformly bounded operators in $H$ for fixet $%
t, $ \ by $\left( 4.1\right) $, $\left( 4.2\right) $ $\left( 4.4\right) $ we
have%
\begin{equation}
\left\Vert \hat{A}^{\frac{1}{2}}BC_{1}\left( t\right) \varphi \right\Vert
_{L^{2}}=\left\Vert F^{-1}\left[ \left\vert \xi \right\vert ^{-1}\hat{g}^{-%
\frac{-1}{2}}\left( \xi \right) \hat{A}^{\frac{1}{2}}C\left( \xi ,t\right) 
\hat{\varphi}\right] \right\Vert _{L^{2}}\lesssim  \tag{4.5}
\end{equation}%
\[
\left\Vert \varphi \right\Vert _{H^{s}\left( A^{\frac{1}{2}}\right) }<\infty
, 
\]%
\[
\left\Vert \hat{A}^{\frac{1}{2}}BS_{1}\left( t\right) \varphi \right\Vert
_{L^{2}}=\left\Vert F^{-1}\left[ \left\vert \xi \right\vert ^{-1}\hat{g}^{-%
\frac{-1}{2}}\left( \xi \right) \hat{A}^{\frac{1}{2}}S\left( \xi ,t\right) 
\hat{\psi}\right] \right\Vert _{L^{2}}\lesssim 
\]%
\[
\left\Vert \psi \right\Vert _{H^{s}\left( A^{\frac{1}{2}}\right) }<\infty . 
\]%
By differentiating $\left( 2.3\right) ,$ in a similar way we have 
\[
\left\Vert \hat{A}^{\frac{1}{2}}B\frac{d}{dt}C_{1}\left( t\right) \varphi
\right\Vert _{L^{2}}=\left\Vert F^{-1}\left[ \left\vert \xi \right\vert ^{-1}%
\hat{g}^{-\frac{-1}{2}}\left( \xi \right) \hat{A}^{\frac{1}{2}}\frac{d}{dt}%
C\left( \xi ,t\right) \hat{\varphi}\right] \right\Vert _{L^{2}} 
\]%
\begin{equation}
\lesssim \left\Vert \varphi \right\Vert _{H^{s}\left( A\right) }<\infty , 
\tag{4.6}
\end{equation}%
\[
\left\Vert \hat{A}^{\frac{1}{2}}B\frac{d}{dt}S_{1}\left( t\right) \varphi
\right\Vert _{L^{2}}=\left\Vert F^{-1}\left[ \left\vert \xi \right\vert ^{-1}%
\hat{g}^{-\frac{-1}{2}}\left( \xi \right) \hat{A}^{\frac{1}{2}}\frac{d}{dt}%
S\left( \xi ,t\right) \hat{\psi}\right] \right\Vert _{L^{2}}\lesssim 
\]%
\[
\left\Vert \psi \right\Vert _{H^{s}\left( A^{\frac{1}{2}}\right) }<\infty . 
\]%
For fixed $t$, we have $f(u)\in H^{s}.$ Moreover,\ \ by assumption on $\hat{A%
}\left( \xi \right) $ we have the uniformly estimate 
\[
\left\Vert \hat{A}^{\frac{1}{2}}\left( \xi \right) \eta ^{-1}\left( \xi
\right) \right\Vert _{B\left( H\right) }\leq C_{A}. 
\]%
\ Then by hypothesis on $\hat{g}\left( \xi \right) ,$ due to $s+r\geq 1$\
from $\left( 4.1\right) $ and $\left( 4.3\right) $\ we get 
\[
\left\Vert \hat{A}^{\frac{1}{2}}BQg\right\Vert _{L^{2}}\leq \left\Vert F^{-1}%
\left[ \left\vert \xi \right\vert \hat{g}^{\frac{1}{2}}\left( \xi \right) 
\hat{A}^{\frac{1}{2}}\left( \xi \right) \dint\limits_{0}^{t}S\left( \xi
,t-\tau \right) \hat{f}\left( u\right) \left( \xi \right) d\tau \right]
\right\Vert _{L^{2}}\lesssim 
\]%
\begin{equation}
C_{A}\left\Vert f\left( u\right) \right\Vert _{H^{s}}<\infty .  \tag{4.7}
\end{equation}%
Then from $\left( 4.2\right) $ and $\left( 4.4\right) -\left( 4.7\right) $
we obtain the assertion.

\textbf{Lemma 4.2.} Assume the Condition 3.2 holds with $a=0$. \ Moreover,
let%
\[
\left\Vert \left( \hat{g}\left( \xi \right) \right) ^{-\frac{-1}{2}%
}\right\Vert _{B\left( H\right) }=O\left( 1+\left\vert \xi \right\vert
^{2}\right) ^{\frac{s+1}{2}}. 
\]%
Suppose the solution of $(1.1)-(1.2)$ exists in $C^{\left( 2\right) }\left( %
\left[ 0,\right. \left. T\right) ;Y_{\infty }^{s,2}\left( A;H\right) \right)
.$ If $B\psi \in L^{2}$ then $Bu_{t}\in C^{\left( 1\right) }\left( \left[
0,\right. \left. T\right) ;L^{2}\right) .$ Moreover, if $B\varphi \in L^{2},$
then $Bu\in C^{\left( 1\right) }\left( \left[ 0,\right. \left. T\right)
;L^{2}\right) .$

\textbf{Proof. }\ Integrating the equation $\left( 1.1\right) $ for $a=0,$
twice and calculating the resulting double integral as an iterated integral,
we have 
\begin{equation}
u\left( x,t\right) =\varphi \left( x\right) +t\psi \left( x\right) - 
\tag{4.8}
\end{equation}%
\[
\dint\limits_{0}^{t}\left( t-\tau \right) \left( A\ast u\right) \left(
x,\tau \right) d\tau +\dint\limits_{0}^{t}\left( t-\tau \right) \Delta
\left( g\ast f\left( u\right) \right) \left( x,\tau \right) d\tau , 
\]

\begin{equation}
u_{t}\left( x,t\right) =\psi \left( x\right) -\dint\limits_{0}^{t}\left(
A\ast u\right) \left( x,\tau \right) d\tau +\dint\limits_{0}^{t}\Delta
\left( g\ast f\left( u\right) \right) \left( x,\tau \right) d\tau . 
\tag{4.9}
\end{equation}%
From $\left( 4.1\right) $ and $\left( 4.9\right) $ for fixed $t$ and $\tau $
we get 
\begin{equation}
\left\Vert Bu_{t}\left( x,t\right) \right\Vert _{L^{2}}=\left\Vert B\psi
\left( x\right) \right\Vert _{L^{2}}-  \tag{4.10}
\end{equation}%
\[
\dint\limits_{0}^{t}\left\Vert B\left( A\ast u\right) \left( x,\tau \right)
\right\Vert _{L^{2}}d\tau -\dint\limits_{0}^{t}\left\Vert B\Delta \left(
g\ast f\left( u\right) \right) \left( x,\tau \right) \right\Vert
_{L^{2}}d\tau . 
\]%
By assumption on $A$, $g$ and by $\left( 4.1\right) $ for fixed $\tau $ we
have%
\[
\left\Vert B\left( A\ast u\right) \left( x,\tau \right) \right\Vert
_{L^{2}}\lesssim \left\Vert F^{-1}\left[ \left\vert \xi \right\vert ^{-1}%
\hat{A}\left( \xi \right) \left( \hat{g}^{-\frac{-1}{2}}\left( \xi \right)
\right) \hat{u}\left( \xi ,\tau \right) \right] \right\Vert _{L^{2}} 
\]%
\begin{equation}
\lesssim \left\Vert u\left( .,\tau \right) \right\Vert _{H^{s}\left(
A\right) }.  \tag{4.11}
\end{equation}%
Moreover, by Lemma 3.3 for all $t$ we have $f\left( u\right) \in H^{s}.$
Also 
\begin{equation}
\left\Vert B\Delta \left( g\ast f\left( u\right) \right) \left( x,\tau
\right) \right\Vert _{L^{2}}\lesssim \left\Vert F^{-1}\left[ \left\vert \xi
\right\vert \left( \hat{g}^{\frac{-1}{2}}\left( \xi \right) \right) \hat{f}%
\left( u\right) \left( \xi \right) \right] \right\Vert _{L^{2}}.  \tag{4.12}
\end{equation}%
Then from $\left( 4.10\right) -\left( 4.12\right) $ we obtain $Bu_{t}\in
C^{\left( 1\right) }\left( \left[ 0,\right. \left. T\right) ;L^{2}\right) .$
The second statement follows similarly from $\left( 4.8\right) .$

From Lemma 4.2 we obtain the following result.

\textbf{Result 4.1.} Assume the Condition 4.1 are satisfied with $a=0$ and%
\[
\left\Vert \hat{g}\left( \xi \right) \right\Vert _{B\left( H\right)
}=O\left( 1+\left\vert \xi \right\vert ^{2}\right) ^{-\frac{r}{2}}. 
\]%
Suppose the solution of $(1.1)-(1.2)$ exists in $C^{\left( 2\right) }\left( %
\left[ 0,T\right] ;Y_{\infty }^{s,2}\left( A;H\right) \right) $ for some $%
s\geq 0.$ If $B\psi \in L^{2}$ then 
\[
Bu_{t}\in C^{\left( 1\right) }\left( \left[ 0,\right. \left. T\right)
;L^{2}\right) . 
\]
Moreover, if $B\varphi \in L^{2},$ then $Bu\in C^{\left( 1\right) }\left( %
\left[ 0,\right. \left. T\right) ;L^{2}\right) .$

Here, 
\begin{equation}
G\left( \sigma \right) =\dint\limits_{0}^{\sigma }f\left( s\right) ds. 
\tag{4.13}
\end{equation}%
\textbf{Lemma 4.3. }Assume the Condition 4.1 is satisfied for $s+r\geq 1$.
Let $B\varphi ,$ $B\psi \in L^{2}$ and $G\left( \varphi \right) \in
L^{1}\left( R^{n};H\right) $. Then for any $t\in \left[ 0,\right. \left.
T\right) $ the energy 
\begin{equation}
E\left( t\right) =\left\Vert Bu_{t}\right\Vert _{L^{2}}^{2}+a\left\Vert
g\ast u\right\Vert _{L^{2}}^{2}+\left\Vert B\left( A\ast u\right)
\right\Vert _{L^{2}}^{2}+2\dint\limits_{R^{n}}G\left( u\right) dx  \tag{4.14}
\end{equation}%
is constant.

\textbf{Proof. } By Theorem 4.1, $\ A^{\frac{1}{2}}Bu,$ $A^{\frac{1}{2}%
}Bu_{t}\in L^{2}.$ By assumptions $g\ast u\in L^{2}$ and $A\ast u\in L^{2}.$
By use of equation $\left( 1.1\right) $, it follows from straightforward
calculation that 
\[
\frac{d}{dt}E\left( t\right) =2\left( Bu_{tt},Bu_{t}\right) +2a\left( F^{-1}%
\hat{g}\ast u,\left( F^{-1}\hat{g}\ast u\right) _{t}\right) + 
\]

\[
2\left[ B\left( A\ast u\right) ,B\left( A\ast u\right) u_{t}\left( t\right) %
\right] +2\left( f\left( u\right) ,u_{t}\right) = 
\]%
\[
2B^{2}\left[ \left( u_{tt}-a\Delta u+A\ast u+\Delta \left[ g\ast f\left(
u\right) \right] ,u_{t}\right) \right] =0, 
\]%
where $\left( u,\upsilon \right) $ denotes the inner product in $L^{2}$
space. \ Integrating the above equality with respect to $t$, we have $\left(
4.14\right) $. By using the above lemmas we obtain the following results

\textbf{Theorem 4.1. }Let the Condition 4.1 holds for $r>2+\frac{n}{2}$. \
Moreover, let $\ B\varphi ,$ $B\psi \in L^{2}$ and $G\left( \varphi \right)
\in L^{1}\left( R^{n};H\right) $ and there is some $k>0$ so that \ $G\left(
\sigma \right) \geq -k\sigma ^{2}$ for all $\sigma \in \mathbb{R}$. Then
there is some $T>0$ such that problem $(1.1)-(1.2)$ has a global solution%
\[
u\in C^{\left( 2\right) }\left( \left[ 0,\right. \left. \infty \right)
;Y_{\infty }^{s,2}\left( A;H\right) \right) . 
\]

\textbf{Proof. }Since $r>2+\frac{n}{2},$ by Theorem 3.4 we get local
existence in 
\[
C^{\left( 2\right) }\left( [0,T);Y_{\infty }^{s,2}\left( A;H\right) \right) 
\]
for some $T>0.$ Assume that $u$ exists on $[0,T).$ By assumption $G\left(
\sigma \right) \geq -k\sigma ^{2}$ and by Lemma 3.4,\ for all $t\in \lbrack
0,T)$ we obtain%
\begin{equation}
\left\Vert Bu_{t}\right\Vert ^{2}+a\left\Vert F^{-1}\hat{g}\ast u\right\Vert
^{2}+\left\Vert B\left( A\ast u\right) \right\Vert ^{2}\leq E\left( 0\right)
+2k\left\Vert u\left( t\right) \right\Vert ^{2}.  \tag{4.15}
\end{equation}

\bigskip By condition on operator function $\hat{g}\left( \xi \right) ,$ we
have 
\begin{equation}
\left\Vert Bu_{t}\right\Vert _{L^{2}\left( A\right)
}^{2}=\dint\limits_{R^{n}}\left\vert \xi \right\vert ^{-2}\left\Vert \hat{g}%
^{-1}\left( \xi \right) \right\Vert _{B\left( H\right) }^{2}\left\Vert A\hat{%
u}_{t}\left( \xi ,t\right) \right\Vert _{H}^{2}\geq  \tag{4.16}
\end{equation}%
\[
C_{g}^{-1}\dint\limits_{R^{n}}\left( 1+\left\vert \xi \right\vert
^{2}\right) ^{\frac{r}{2}-1}\left\Vert A\hat{u}_{t}\left( \xi ,t\right)
\right\Vert _{H}^{2}\approx C_{g}^{-1}\left\Vert Au_{t}\left( t\right)
\right\Vert _{H^{\frac{r}{2}-1}}^{2}, 
\]%
$C_{g}$ is the positive constant that appears in $\left( 3.19.\right) $. By
properties of norms in Hilbert spaces and by Cauchy-Schwarz inequality, from 
$\left( 4.15\right) $ and $\left( 4.16\right) $ we get 
\[
\frac{d}{dt}\left\Vert u\left( t\right) \right\Vert _{H^{\frac{r}{2}%
-1}\left( A\right) }^{2}\leq 2\left\Vert u_{t}\left( t\right) \right\Vert
_{H^{\frac{r}{2}-1}\left( A\right) }\left\Vert u\left( t\right) \right\Vert
_{H^{\frac{r}{2}-1}\left( A\right) }\leq 
\]%
\[
\left\Vert u_{t}\left( t\right) \right\Vert _{H^{\frac{r}{2}-1}\left(
A\right) }^{2}+\left\Vert u\left( t\right) \right\Vert _{H^{\frac{r}{2}%
-1}\left( A\right) }^{2}\leq C\left\Vert Bu_{t}\left( t\right) \right\Vert
_{H^{\frac{r}{2}-1}\left( A\right) }^{2}+ 
\]%
\[
\left\Vert u\left( t\right) \right\Vert _{H^{\frac{r}{2}-1}\left( A\right)
}^{2}\leq CE\left( 0\right) +\left( 2Ck+1\right) \left\Vert u\left( t\right)
\right\Vert _{H^{\frac{r}{2}-1}\left( A\right) }^{2} 
\]

Gronwall's lemma implies that $\left\Vert u\left( t\right) \right\Vert _{H^{%
\frac{r}{2}-1}\left( A\right) }$ is bounded in $[0,T)$. But, since $(r/2)-1>%
\frac{n}{4}$, we conclude that $\left\Vert u(t)\right\Vert _{L^{\infty
}\left( A\right) }$ also is bounded in $[0,T)$. By Lemma $3.8$ this implies
a global solution.

\begin{center}
\textbf{\ \ 5. Applications}

\textbf{5.1}.\textbf{The Cauchy problem for the system of nonlocal WEs }
\end{center}

\bigskip Consider the problem $\left( 1.3\right) $. Let%
\[
\text{ }l_{2}\left( N\right) =\left\{ \text{ }u=\left\{ u_{j}\right\} ,\text{
}j=1,2,...N,\left\Vert u\right\Vert _{l_{2}\left( N\right) }=\left(
\sum\limits_{j=1}^{N}\left\vert u_{j}\right\vert ^{2}\right) ^{\frac{1}{2}%
}<\infty \right\} , 
\]%
where $N\in \mathbb{N}$ (see $\left[ \text{23, \S\ 1.18}\right] .$ Let $%
A_{1} $ be the operator in $l_{2}\left( N\right) $ defined by%
\[
\text{ }A_{1}=\left[ a_{jm}\left( x\right) \right] ,\text{ }%
a_{jm}=b_{j}\left( x\right) 2^{\sigma m},\text{ }m,j=1,2,...,N,\text{ }%
D\left( A_{1}\right) =\text{ }l_{2}^{\sigma }\left( N\right) = 
\]

\[
\left\{ \text{ }u=\left\{ u_{j}\right\} ,\text{ }j=1,2,...N,\left\Vert
u\right\Vert _{l_{2}^{\sigma }\left( N\right) }=\left(
\sum\limits_{j=1}^{N}2^{\sigma j}u_{j}^{2}\right) ^{\frac{1}{2}}<\infty
\right\} ,\text{ }\sigma >0. 
\]

Let \ 
\[
H^{s,p}\left( E\right) =H^{s,p}\left( R^{n};E\right) ,\text{ }H^{s}\left(
E\right) =H^{s,2}\left( R^{n};E\right) , 
\]%
\[
Y^{s,p,\sigma }=H^{s,p}\left( R^{n};l_{2}\right) \cap L^{p}\left(
R^{n};l_{2}^{\sigma }\right) ,\text{ }1\leq q\leq \infty , 
\]%
\[
H_{0}\left( l_{2}\right) =H^{s\left( 1-\frac{1}{2p}\right) ,p}\left(
R^{n};l_{2}\right) \cap L^{p}\left( R^{n};l_{2}^{\sigma \left( 1-\frac{1}{2p}%
\right) }\right) . 
\]

\ Let $f=\left\{ f_{m}\right\} $, $m=1,2,...,N$ and let%
\[
\eta _{1}=\eta _{1}\left( \xi \right) =\left[ \hat{a}\left( \xi \right)
\left\vert \xi \right\vert ^{2}+\hat{A}_{1}\left( \xi \right) \right] ^{%
\frac{1}{2}}.
\]%
\ From Theorem 3.1 we obtain the following result \ 

\textbf{Theorem 5.1.} Assume: (1) $a\in L^{1}\left( R^{n}\right) ,$ $\hat{a}%
\left( \xi \right) >0$ for all $\xi \in R^{n}$, $\hat{b}_{j}=b_{j}\left( \xi
\right) $ are nonnegat\i ve bounded d\i fferentiable functions on $R^{n}$
and $\hat{a}\left( \xi \right) +\hat{b}_{j}\left( \xi \right) >0$ for $\xi
\in R^{n};$ (2) $D^{\alpha }\hat{b}_{j}$ are uniformly bounded on $R^{n}$
for $\alpha =\left( \alpha _{1},\alpha _{2},...,\alpha _{n}\right) ,$ $%
\left\vert \alpha \right\vert \leq n$ and the uniform estimate holds 
\[
\dsum\limits_{j=1}^{N}\left\vert D^{\alpha }\hat{b}_{j}\left( \xi \right)
\right\vert ^{2}\left[ \hat{a}\left( \xi \right) +\hat{b}_{j}\left( \xi
\right) \right] ^{-1}\leq M;
\]

\bigskip (3) $\varphi \in H^{s,p}\left( l_{2}^{\sigma }\right) ,$ $\psi \in
H^{s,p}\left( l_{2}^{\frac{\sigma }{2}}\right) $\ and $s>1+\frac{n}{p}$ for $%
p\in \left[ 1,\infty \right] $; (4) the kernel $g_{mj}$ are bounded
integrable functions, whose Fourier transform satisfies%
\[
0\leq \dsum\limits_{j=m,j}^{N}\left\vert \hat{g}_{mj}\left( \xi \right)
\right\vert ^{2}\lesssim \left( 1+\left\vert \xi \right\vert ^{2}\right) ^{-%
\frac{r}{2}}\text{ for all }\xi \in R^{n}\text{ and }r\geq 2. 
\]%
(4) the function%
\[
u\rightarrow f\left( x,t,u\right) :R^{n}\times \left[ 0,T\right] \times
H_{0}\left( l_{2}\right) \rightarrow l_{2} 
\]%
is a measurable in $\left( x,t\right) \in R^{n}\times \left[ 0,T\right] $
for $u\in H_{0}\left( l_{2}\right) ;$ Moreover, $f\left( x,t,u\right) $ is
continuous in $u\in H_{0}\left( l_{2}\right) $ and $f\left( x,t,u\right) \in
C^{\left[ s\right] +1}\left( H_{0}\left( l_{2}\right) ;l_{2}\right) $
uniformly in $x\in R^{n},$ $t\in \left[ 0,T\right] .$

\bigskip Then problem $\left( 1.3\right) $ has a unique local strange
solution%
\[
u\in C^{\left( 2\right) }\left( \left[ 0,\right. \left. T_{0}\right)
;Y_{\infty }^{s,p}\left( A_{1},H\right) \right) , 
\]%
where $T_{0}$ is a maximal time interval that is appropriately small
relative to $M$. Moreover, if

\[
\sup_{t\in \left[ 0\right. ,\left. T_{0}\right) }\left( \left\Vert
u\right\Vert _{Y_{\infty }^{s,p}\left( A_{1};H\right) }+\left\Vert
u_{t}\right\Vert _{Y_{\infty }^{s,p}\left( A_{1};H\right) }\right) <\infty 
\]%
then $T_{0}=\infty .$

\ \textbf{Proof. }By virtue of $\left[ \text{23, \S\ 1.18}\right] ,$ $%
l_{2}\left( N\right) $ is a Hilbert space. By definition of $H^{s,p}\left(
A_{1},l_{2}\right) $ and by real interpolation of Banach spaces (see e.g. $%
\left[ \text{23, \S 1.3, 1.18}\right] $) we have%
\[
H_{0}\left( l_{2}\right) =\left( Y^{s,p,\sigma },L^{p}\left(
R^{n};l_{2}\right) \right) _{\frac{1}{2p},p}= 
\]%
\[
\left( H^{s,p}\left( R^{n};l_{2}\right) \cap L^{p}\left( R^{n};l_{2}^{\sigma
}\right) ,L^{p}\left( R^{n};l_{2}\right) \right) _{\frac{1}{2p},p}= 
\]

\[
\left( H^{s,p}\left( R^{n};l_{2}\right) ,L^{p}\left( R^{n};l_{2}\right)
\right) _{\frac{1}{sp},p}\cap \left( L^{p}\left( R^{n};l_{2}^{\sigma
}\right) ,L^{p}\left( R^{n};l_{2}\right) \right) _{\frac{1}{2p},p}= 
\]%
\[
H^{s\left( 1-\frac{1}{2p}\right) ,p}\left( R^{n};l_{2}\right) \cap
L^{p}\left( R^{n};l_{2}^{\sigma \left( 1-\frac{1}{sp}\right) }\right)
=H_{0}\left( l_{2}\right) . 
\]%
By assumptions (1), (2) we obtain that $\hat{A}_{1}\left( \xi \right) $ is
uniformly sectorial in $l_{2}$, $\hat{A}_{1}\left( \xi \right) \in \sigma
\left( M_{0},\omega ,l_{2}\right) ,$ $\eta _{1}\left( \xi \right) \neq 0$\
for all $\ \xi \in R^{n}$\ and 
\[
\left\Vert D^{\alpha }\hat{A}_{1}\left( \xi \right) \eta _{1}^{-1}\left( \xi
\right) \right\Vert _{B\left( l_{2}\right) }\leq M 
\]%
for $\alpha =\left( \alpha _{1},\alpha _{2},...,\alpha _{n}\right) ,$ $%
\left\vert \alpha \right\vert \leq n.$\ Hence, by (3), (4) all conditions of
Theorem 3.1 are hold, i,e., we get the conclusion.\ 

Let $G$ be a function defined by $\left( 4.15\right) .$

\textbf{Theorem 5.2.} Assume: (a) (1)-(3) assumptions of Theorem 5.1 are
satisfied for $p=2$; (b) $\varphi \in H^{s}\left( l_{2}^{\sigma }\right) ,$ $%
\psi $ $\in H^{s}\left( l_{2}^{\frac{\sigma }{2}}\right) $ for $s>1+\frac{n}{%
2}$ and 
\[
\left\Vert \hat{g}\left( \xi \right) \right\Vert _{B\left( l_{2}\right)
}\lesssim \left( 1+\left\vert \xi \right\vert ^{2}\right) ^{-\frac{r}{2}}%
\text{ for }r\leq 2\left( s+1\right) , 
\]%
\[
\left\Vert \hat{g}^{\frac{1}{2}}\left( \xi \right) \right\Vert _{B\left(
l_{2}\right) }\lesssim \left\vert \xi \right\vert \left( 1+\left\vert \xi
\right\vert ^{2}\right) ^{\frac{s}{2}}\text{ for all }\xi \in R^{n}; 
\]%
(c) $f_{m}\in C^{\left[ s\right] }\left( R;l_{2}\right) $ with $f(0)=0$ and%
\[
\dsum\limits_{m=1}^{N}\left\vert \hat{f}_{m}\left( u\right) \left( \xi
\right) \right\vert ^{2}<\infty \text{ for all }u=\left(
u_{1},u_{2},....,u_{m}\right) \in C^{\left( 2\right) }\left( \left[
0,\right. \left. \infty \right) ;Y_{\infty }^{s,2}\left( A_{2};l_{2}\right)
\right) ; 
\]

(d) $B\varphi $, $B\psi \in L^{2}\left( R^{n};l_{2}\right) $ and $G\left(
\varphi \right) \in L^{1}\left( R^{n};l_{2}\right) $; (e) there is some $k>0$
so that $G\left( \nu \right) \geq -k\nu ^{2}$ for all $\nu \in \mathbb{R}$.
Then there is some $T>0$ such that problem $(1.3)$ has a global solution%
\[
u\in C^{\left( 2\right) }\left( \left[ 0,\right. \left. \infty \right)
;Y_{\infty }^{s,2}\left( A_{1};l_{2}\right) \right) . 
\]

\ \textbf{Proof. }From the assumptions (a), (b) it is clear to see that the
Condition 4.1 holds for $H=l_{2}$ and $r>2+\frac{n}{2}.$ By (c), (d), (e)
all other assumptions of Theorem 4.1 are satisfied. Hence, we obtain the
assertion.

\begin{center}
\textbf{5.2.} \textbf{The mixed problem for degenerate nonlocal WE}
\end{center}

Consider the problem $\left( 1.5\right) -\left( 1.7\right) $. \ Let%
\[
Y^{s,p,2}=H^{s,p}\left( R^{n};L^{2}\left( 0,1\right) \right) \cap
L^{p}\left( R^{n};H^{\left[ 2\right] }\left( 0,1\right) \right) ,\text{ }%
1\leq p\leq \infty , 
\]%
\[
\text{ }H_{0}\left( L^{2}\right) =H^{s\left( 1-\frac{1}{2p}\right) }\left(
R^{n};L^{2}\left( 0,1\right) \right) \cap L^{p}\left( R^{n};H^{\left[
2\left( 1-\frac{1}{2p}\right) \right] }\left( 0,1\right) \right) . 
\]

Let $A_{2}$ is the operator in $L^{2}\left( 0,1\right) $ defined by $\left(
1.4\right) $ and let%
\[
\eta _{2}=\eta _{2}\left( \xi \right) =\left[ \hat{a}\left( \xi \right)
\left\vert \xi \right\vert ^{2}+\hat{A}_{2}\left( \xi \right) \right] ^{%
\frac{1}{2}}.
\]%
Now, we present the following result:

\textbf{Condition 5.1 }Assume;

(1) $0\leq \gamma <\frac{1}{2},$ $\alpha _{10}\beta _{20}-\alpha _{20}\beta
_{10}\neq 0,$ $\alpha _{20}\beta _{11}+\alpha _{21}\beta _{10}-\alpha
_{10}\beta _{21}-\alpha _{11}\beta _{20}\neq 0,$ $\alpha _{11}\beta
_{21}+\alpha _{21}\beta _{11}\neq 0$, $\alpha _{11}\beta _{21}-\alpha
_{11}\alpha _{21}\neq 0,$ $\alpha _{11}\neq \beta _{11}$ for $\nu _{k}=1$
and $\left\vert \alpha _{k0}\right\vert +\left\vert \beta _{k0}\right\vert
>0,$ $\alpha _{10}\alpha _{20}+\beta _{10}\beta _{20}\neq 0$ for $\nu
_{k}=0; $

(2) $b_{1}$ and $b_{2}$ are complex valued functions on $\left( 0,1\right) $%
. Morover, $b_{1}\in C\left[ 0,1\right] ,$ $b_{1}\left( 0\right)
=b_{1}\left( 1\right) $, $b_{2}\in L_{\infty }\left( 0,1\right) $ and $%
\left\vert b_{2}\left( x\right) \right\vert \leq C$ $\left\vert b_{1}^{\frac{%
1}{2}-\mu }\left( x\right) \right\vert $ for $0<\mu <\frac{1}{2}$ and for
a.a. $x\in \left( 0,1\right) ;$

(3) $a\geq 0,$ $D^{\alpha }\hat{b}_{j},$ $j=1,$ $2$ are uniformly bounded on 
$R^{n}$ for all $\alpha =\left( \alpha _{1},\alpha _{2},...,\alpha
_{n}\right) $ with $\left\vert \alpha \right\vert \leq n$ and $\eta
_{2}\left( \xi \right) \neq 0$\ for all $\xi \in R^{n};$

(4) $\varphi \in W^{s,p}\left( A_{2}\right) $ and $\psi \in W^{s,p}\left(
A_{2}^{\frac{1}{2}}\right) ;$

(5)$\ $for $\alpha =\left( \alpha _{1},\alpha _{2},...,\alpha _{n}\right) ,$ 
$\left\vert \alpha \right\vert \leq n$ the uniform estimate holds 
\[
\left\Vert \left[ D^{\alpha }\hat{A}_{2}\left( \xi \right) \right] \eta
_{2}^{-1}\left( \xi \right) \right\Vert _{B\left( H\right) }\leq M. 
\]

(6) the function%
\[
u\rightarrow f\left( x,t,u\right) :R^{n}\times \left[ 0,T\right] \times
H_{0}\left( L^{p_{1}}\left( 0,1\right) \right) \rightarrow L^{2}\left(
0,1\right) 
\]%
is a measurable in $\left( x,t\right) \in R^{n}\times \left[ 0,T\right] $
for $u\in H_{0}\left( L^{2}\left( 0,1\right) \right) ;$ $f\left(
x,t,u\right) $. Moreover, $f\left( x,t,u\right) $ is continuous in $u\in
H_{0}\left( L^{2}\left( 0,1\right) \right) $ and%
\[
f\left( x,t,u\right) \in C^{\left[ s\right] +1}\left( H_{0}\left(
L^{2}\left( 0,1\right) \right) ;L^{2}\left( 0,1\right) \right) 
\]%
uniformly with respect to $x\in R^{n},$ $t\in \left[ 0,T\right] .$

\bigskip \textbf{Theorem 5.3.} Assume that the Condition 5.1 is satisfied.
Suppose $\varphi \in W^{s,p}\left( R^{n};W^{2,2}\left( 0,1\right) \right) ,$ 
$\psi \in W^{s,p}\left( R^{n};W^{1,2}\left( 0,1\right) \right) $\ for $s>1+%
\frac{n}{p}$ and $p\in \left[ 1,\infty \right] $.

Then problem $\left( 1.5\right) -\left( 1.7\right) $ has a unique local
strange solution 
\[
u\in C^{\left( 2\right) }\left( \left[ 0,\right. \left. T_{0}\right)
;Y_{\infty }^{s,p}\left( A_{2},L^{2}\left( 0,1\right) \right) \right) , 
\]%
where $T_{0}$ is a maximal time interval that is appropriately small
relative to $M$. Moreover, if

\[
\sup_{t\in \left[ 0\right. ,\left. T_{0}\right) }\left( \left\Vert
u\right\Vert _{Y_{\infty }^{s,p}\left( A_{2};L^{2}\left( 0,1\right) \right)
}+\left\Vert u_{t}\right\Vert _{Y_{\infty }^{s,p}\left( A_{2};L^{2}\left(
0,1\right) \right) }\right) <\infty 
\]%
then $T_{0}=\infty .$

\ \textbf{Proof.}\ It is known that $L^{2}\left( 0,1\right) $ is a Hilbert
space. By definition of $H^{s,p}\left( A_{2},L^{2}\left( 0,1\right) \right) $
and by real interpolation of Banach spaces (see e.g. $\left[ \text{23, \S 1.3%
}\right] $) we have%
\[
\left( H^{s,p}\left( A_{2},L^{2}\left( 0,1\right) \right) ,L^{p}\left(
R^{n};L^{2}\left( 0,1\right) \right) \right) _{\frac{1}{2p},p}= 
\]%
\[
\left( H^{s,p}\left( R^{n};L^{2}\left( 0,1\right) \right) \cap L^{p}\left(
R^{n};H^{\left[ 2\right] ,2}\left( 0,1\right) \right) ,L^{p}\left(
R^{n};L^{2}\left( 0,1\right) \right) \right) _{\frac{1}{2p},p}= 
\]%
\[
\left( H^{s,p}\left( R^{n};L^{2}\left( 0,1\right) \right) ,L^{p}\left(
R^{n};L^{2}\left( 0,1\right) \right) \right) _{\frac{1}{2p},p}\cap 
\]%
\[
\left( L^{p}\left( R^{n};H^{\left[ 2\right] ,2}\left( 0,1\right) \right)
,L^{p}\left( R^{n};L^{2}\left( 0,1\right) \right) \right) _{\frac{1}{2p},p}= 
\]%
\[
H^{s\left( 1-\frac{1}{2p}\right) }\left( R^{n};L^{2}\left( 0,1\right)
\right) \cap L^{p}\left( R^{n};H^{\left[ 2\left( 1-\frac{1}{2p}\right) %
\right] }\left( 0,1\right) \right) =H_{0}\left( L^{2}\right) . 
\]%
In view of $\left[ \text{32, Theorem 4.1}\right] $ we obtain that $\hat{A}%
_{2}\left( \xi \right) $ is uniformly sectorial in $L^{2}\left( 0,1\right) $
and 
\[
\hat{A}_{2}\left( \xi \right) \in \sigma \left( M_{0},\omega ,L^{2}\left(
0,1\right) \right) . 
\]

Moreover, by using the assumptions (1), (2) we deduced that $\eta _{2}\left(
\xi \right) \neq 0$\ for all $\ \xi \in R^{n}$\ and 
\[
\left\Vert D^{\alpha }\hat{A}_{2}\left( \xi \right) \eta _{2}^{-1}\left( \xi
\right) \right\Vert _{B\left( l_{2}\right) }\leq M. 
\]%
for $\alpha =\left( \alpha _{1},\alpha _{2},...,\alpha _{n}\right) ,$ $%
\left\vert \alpha \right\vert \leq n.$ Hence, by hypothesis (3), (4) of the
Condition 5.1 we get that all hypothesis of Theorem 4.1 are hold, i,e., we
obtain the conclusion.\ 

\bigskip \textbf{Theorem 5.4.} Assume the Condition 5.1 is satisfied.
Suppose 
\[
\varphi \in H^{s,2}\left( R^{n};H^{\left[ 2\right] ,2}\left( 0,1\right)
\right) ,\psi \in H^{s,2}\left( R^{n};H^{\left[ 1\right] ,2}\left(
0,1\right) \right) \text{, }s>1+\frac{n}{2}. 
\]%
\ Suppose $f\in C^{\left[ s\right] }\left( R;L^{2}\left( \left( 0,1\right)
\right) \right) $ with $f(0)=0.$\ Let the kernel $g_{mj}$ be bounded
integrable functions and

\[
\left\Vert \hat{g}\left( \xi \right) \right\Vert _{B\left( l_{2}\right)
}\lesssim \left( 1+\left\vert \xi \right\vert ^{2}\right) ^{-\frac{r}{2}}%
\text{ for }r\leq 2\left( s+1\right) , 
\]%
\[
\left\Vert \hat{g}^{\frac{1}{2}}\left( \xi \right) \right\Vert _{B\left(
l_{2}\right) }\lesssim \left\vert \xi \right\vert \left( 1+\left\vert \xi
\right\vert ^{2}\right) ^{\frac{s}{2}}\text{ for all }\xi \in R^{n} 
\]

Moreover, let $\ B\varphi ,$ $B\psi \in L^{2}\left( \left( 0,1\right) \times
R^{n}\right) $, $G\left( \varphi \right) \in L^{1}\left( \left( 0,1\right)
\times R^{n}\right) $ and there is some $k>0$ so that \ $G\left( r\right)
\geq -kr^{2}$ for $r\in \mathbb{R}$. Then there is some $T>0$ such that the
problem $(1.5)-\left( 1.7\right) $ has a global solution%
\[
u\in C^{2}\left( \left[ 0,\right. \left. \infty \right) ;Y_{\infty
}^{s,2}\right) . 
\]

\textbf{Proof.} Indeed, by assumptions all conditions of Theorem 4.1. are
satisfied for $H=L^{2}\left( 0,1\right) ,$ i.e. we obtain the assertion.

\textbf{References}

\begin{quote}
\ \ \ \ \ \ \ \ \ \ \ \ \ \ \ \ \ \ \ \ \ \ \ \ 
\end{quote}

\begin{enumerate}
\item A. C. Eringen, Nonlocal Continuum Field Theories, New York, Springer
(2002),.

\item Z. Huang, Formulations of nonlocal continuum mechanics based on a new
definition of stress tensor Acta Mech. (2006)187, 11--27.

\item C. Polizzotto, Nonlocal elasticity and related variational principles
Int. J. Solids Struct. ( 2001) 38 7359--80.

\item C. A. Silling, Reformulation of elasticity theory for discontinuities
and long-range forces J. Mech. Phys. Solids (2000)48 175-209.

\item M. Arndt and M. Griebel, Derivation of higher order gradient continuum
models from atomistic models for crystalline solids Multiscale Modeling
Simul. (2005)4, 531--62..

\item X. Blanc, C. LeBris, P. L. Lions, Atomistic to continuum limits for
computational materials science, ESAIM--- Math. Modelling Numer. Anal.
(2007)41, 391--426.

\item A. De Godefroy, Blow up of solutions of a generalized Boussinesq
equation IMA J. Appl. Math.(1998) 60 123--38.

\item A. Constantin and L. Molinet, The initial value problem for a
generalized Boussinesq equation, Diff.Integral Eqns. (2002)15, 1061--72.

\item G. Chen and S. Wang, Existence and nonexistence of global solutions
for the generalized IMBq equation Nonlinear Anal.---Theory Methods Appl.
(1999)36, 961--80.

\item M. Lazar, G. A. Maugin and E. C. Aifantis, On a theory of nonlocal
elasticity of bi-Helmholtz type and some applications Int. J. Solids and
Struct. (2006)43, 1404--21.

\item N. Duruk, H.A. Erbay and A. Erkip, Global existence and blow-up for a
class of nonlocal nonlinear Cauchy problems arising in elasticity,
Nonlinearity, (2010)23, 107--118.

\item S. Wang, G. Chen, Small amplitude solutions of the generalized IMBq
equation, J. Math. Anal. Appl. 274 (2002) 846--866.

\item S.Wang and G.Chen, Cauchy problem of the generalized double dispersion
equation Nonlinear Anal.--- Theory Methods Appl. (2006 )64 159--73.

\item J.L. Bona, R.L. Sachs, Global existence of smooth solutions and
stability of solitary waves for a generalized Boussinesq equation, Comm.
Math. Phys. 118 (1988), 15--29.

\item F. Linares, Global existence of small solutions for a generalized
Boussinesq equation, J. Differential Equations 106 (1993), 257--293.

\item Y. Liu, Instability and blow-up of solutions to a generalized
Boussinesq equation, SIAM J. Math. Anal. 26 (1995), 1527--1546.

\item V.G. Makhankov, Dynamics of classical solutions (in non-integrable
systems), Phys. Lett. C 35(1978), 1--128.

\item G.B. Whitham, Linear and Nonlinear Waves, Wiley--Interscience, New
York, 1975.

\item N.J. Zabusky, Nonlinear Partial Differential Equations, Academic
Press, New York, 1967.

\item A. Ashyralyev, N. Aggez, Nonlocal boundary value hyperbolic problems
involving Integral conditions, Bound.Value Probl., 2014 V (2014):214.

\item L. S. Pulkina, A non local problem with integral conditions for
hyperbolice quations, Electron.J.Differ.Equ.(1999)45, 1-6.

\item M. Girardi, L. Weis, Operator-valued multiplier theorems on Besov
spaces, Math. Nachr. 251 (2003), 34-51.

\item H. Triebel, Interpolation theory, Function spaces, Differential
operators, North-Holland, Amsterdam, 1978.

\item H. Triebel, Fractals and spectra, Birkhauser Verlag, Related to
Fourier analysis and function spaces, Basel, 1997.

\item L. Nirenberg, On elliptic partial differential equations, Ann. Scuola
Norm. Sup. Pisa (1959)13 , 115--162.

\item S. Klainerman, Global existence for nonlinear wave equations, Comm.
Pure Appl. Math.(1980)33 , 43--101.

\item R. Coifman and Y. Meyer, Wavelets. Calder%
\'{}%
on-Zygmund and Multilinear Operators, Cambridge University Press, 1997.

\item T. Kato, G. Ponce, Commutator estimates and the Euler and
Navier--Stokes equations, Comm. Pure Appl. Math. (1988)41, 891--907.

\item H. O. Fattorini, Second order linear differential equations in Banach
spaces, in North Holland Mathematics Studies,\ V. 108, North-Holland,
Amsterdam, 1985.

\item A. Pazy, Semigroups of linear operators and applications to partial
differential equations. Springer, Berlin, 1983.

\item V. B. Shakhmurov, Embedding and separable differential operators in
Sobolev-Lions type spaces, Math. Notes, 84(2008) (6), 906-926.

\item V. B. Shakhmurov, Linear and nonlinear abstract differential equations
with small marameters, Banach J. Math. Anal. 10 (2016)(1), 147--168.
\end{enumerate}

\end{document}